\documentclass[10pt,leqno]{article}

\usepackage[francais]{babel}
\usepackage{amsmath,amssymb,amscd}
\usepackage{a4wide}

\usepackage[latin1]{inputenc}

\usepackage[all]{xy} 

\makeatletter
\@addtoreset{equation}{subsection}         
\makeatother

\newenvironment{paragr}[1][]{\refstepcounter{subsection} \noindent \textbf{\thesubsection . \ #1}}{\medskip}

\newenvironment{theoreme}{ \medskip\refstepcounter{theo}  \noindent\textbf{Th\'eor\`eme \thetheo}. ---\em}{\em \medskip}
\newenvironment{proposition}{\medskip\refstepcounter{theo}   \noindent\textbf{Proposition \thetheo}. ---\em}{\em\medskip}
\newenvironment{corollaire}{\medskip\refstepcounter{theo}  \noindent\textbf{Corollaire \thetheo}. ---\em}{\em\medskip}

\newenvironment{lemme}{\medskip\refstepcounter{theo}   \noindent\textbf{Lemme \thetheo}. ---\em}{\em\medskip}

\newenvironment{preuve}[1][]{\noindent \textbf{Démonstration.} #1 --- }{\hfill
  \ensuremath{\square} \medskip}

\newenvironment{remarque}{\medskip\refstepcounter{theo}  \noindent\textbf{Remarque \thetheo}. ---}{\medskip}
\newenvironment{remarques}{\medskip\refstepcounter{theo}  \noindent\textbf{Remarques \thetheo}. ---}{\medskip}

\DeclareMathOperator{\vol}{vol}

\DeclareMathOperator{\Ad}{Ad}

\DeclareMathOperator{\ad}{ad}

\DeclareMathOperator{\der}{der}

\DeclareMathOperator{\Norm}{Norm}

\DeclareMathOperator{\Hom}{Hom}

\DeclareMathOperator{\Ind}{Ind}

\DeclareMathOperator{\Res}{Res}

\newcommand{\ZZ}{\mathbb{Z}}
\newcommand{\Gm}{\mathbb{G}_m}

\newcommand{\NN}{\mathbb{N}}
\newcommand{\RR}{\mathbb{R}}
\newcommand{\AAA}{\mathbb{A}}
\newcommand{\CC}{\mathbb{C}}

\newcommand{\QQ}{\mathbb{Q}}

\newcommand{\oc}{\mathcal{O}}

\newcommand{\Sc}{\mathcal{S}}

\newcommand{\lc}{\mathcal{L}}
\newcommand{\dc}{\mathcal{D}}
\newcommand{\fc}{\mathcal{F}}
\newcommand{\pc}{\mathcal{P}}

\newcommand{\kc}{\mathcal{K}}

\newcommand{\ggo}{\mathfrak{g}}

\newcommand{\of}{\mathfrak{o}}

\newcommand{\mgo}{\mathfrak{m}}
\newcommand{\ngo}{\mathfrak{n}}

\newcommand{\pgo}{\mathfrak{p}}
\newcommand{\qgo}{\mathfrak{q}}

\newcommand{\hgo}{\mathfrak{h}}

\newcommand{\lgo}{\mathfrak{l}}

\newcommand{\al}{\alpha}
\newcommand{\be}{\beta}

\newcommand{\om}{\omega}

\newcommand{\la}{\lambda}

\newcommand{\back}{\backslash}

\newcommand{\Cc}{C_c^\infty}

\newcommand{\bg}{\langle}
\newcommand{\bd}{\rangle}

\newcommand{\eps}{\varepsilon}

\renewcommand{\leq}{\leqslant}
\renewcommand{\geq}{\geqslant}

\title{Sur certaines contributions unipotentes dans la formule des traces d'Arthur}
\author{Pierre-Henri Chaudouard}
\date{}
\begin{document}
\maketitle

\begin{abstract}
  Nous établissons un développement fin pour la partie géométrique de la formule des traces d'Arthur tel que conjecturé par Werner Hoffmann. Pour le groupe général linéaire, nous en déduisons une expression pour la contribution des orbites unipotentes régulières par blocs (ce sont des orbites avec un seul bloc de Jordan avec éventuellement une multiplicité). Comme conséquence, nous obtenons des formules explicites pour les coefficients globaux d'Arthur attachés à de telles orbites.
\end{abstract}

\renewcommand{\abstractname}{Abstract}
\begin{abstract}
  We  establish a fine  expansion for the geometric part of the Arthur-Selberg trace formula (as it was conjectured by Werner Hoffmann). For the general linear group, we deduce an expression for the contributions of regular by blocks unipotent orbits (orbits with one Jordan block with multiplicity).  As a consequence, we find formulas for Arthur's global coefficients attached to such orbits.
\end{abstract}

\tableofcontents

\section{Introduction}

\begin{paragr}
  La formule des traces d'Arthur-Selberg joue un rôle central en théorie des formes automorphes. 
De nombreuses applications reposent sur une comparaison de deux formules des traces. Néanmoins on peut espérer tirer une information intéressante d'une seule formule des traces (par exemple une dimension d'espaces de formes automorphes). Mais pour cela, il faut disposer d'une formule aussi explicite que possible. Dans l'état actuel de la formule des traces, son côté géométrique est une combinaison de distributions locales avec des coefficients, de nature globale, dont on ne connait pas de formules sauf pour les termes semi-simples (ce sont alors des \og volumes\fg{}). Comme il y a une sorte de décomposition de Jordan pour ces coefficients, le point crucial est de les comprendre pour les termes unipotents. Dans cette article, on résout la question pour le groupe $GL(n)$ et les orbites des éléments unipotents dont la décomposition de Jordan ne comprend qu'un seul bloc avec éventuellement une multiplicité (on les appelle réguliers par blocs).  L'article \cite{scfhn}, dont le présent article généralise les méthodes,  avait abordé un problème similaire pour les corps de fonctions sous un angle un peu différent. Dans l'article  \cite{scuft}, on avait également donné une expression pour les coefficients des  orbites des éléments unipotents dont la décomposition de Jordan comprend \emph{tous} les blocs de taille $1$ à une taille donnée. Signalons les  articles récents \cite{Matz-bounds} et \cite{HW} où l'étude de ces coefficients est également abordée. 

Notre point de départ est un développement fin du côté géométrique de la formule des traces (un tel développement avait été conjecturé par Hoffmann dans \cite{Hoff}). 
\end{paragr}

\begin{paragr}[Développement fin de la formule des traces.] --- Dans cet article, on travaille en fait sur les algèbres de Lie (c'est-à-dire qu'on remplace l'action par conjugaison d'un groupe sur lui-même par son action adjointe). Dans ce cadre, on dispose d'un analogue de la formule des traces d'Arthur (cf. \cite{PH1}) dont le développement nilpotent comprend (via un passage par l'exponentielle) les coefficients globaux unipotents mentionnés précédemment. Soit $G$ un groupe réductif sur un corps de nombres $F$ et $\ggo$ son algèbre de Lie. Soit $\AAA$ l'anneau des adèles de $F$ et $f\in \Sc(\ggo(\AAA))$ une fonction de Bruhat-Schwartz. Soit $\of$ une orbite \og géométrique\fg{} de $G$ agissant sur $\ggo$. Pour chaque sous-groupe parabolique standard $P$ de décomposition de Levi $MN$, on a un ersatz de noyau (une fonction sur $M(F)N(\AAA)\back G(\AAA)$)
$$k_{P,\of}(f,g)=\sum_{X\in \mgo(F)_\of} \int_{\ngo(\AAA)} f(\Ad(g^{-1}(X+U))\,dU
$$
où  $\mgo(F)_\of$ est l'ensemble des $X\in \mgo(F)$ tels que \og l'orbite induite\fg{} soit $\of$. En s'inspirant des constructions d'Arthur, on forme alors une sorte de noyau automorphe modifié (une fonction sur $G(F)\back G(\AAA)$)
$$
k_\of^T(f,g)=\sum_P \eps_P^G \sum_{\delta \in P(F)\back G(F)} \hat{\tau}_P(\delta g,T) k_{P,\of}(\delta g, f).
$$
où la somme porte sur les sous-groupes paraboliques standard. Le coefficient  $\eps_P^G$ est un signe et $\hat{\tau}_P $ est une certaine fonction caractéristique qui assure, entre autres, que la somme sur $\delta$ est à support fini. Enfin, $T$ est un paramètre de troncature.
Voici le résultat général qu'on obtient (on note $G(\AAA)^1$  un certain sous-groupe de $G(\AAA)$ pour lequel le quotient $G(F)\back G(\AAA)^1$ est de volume fini).

\begin{theoreme} (pour des énoncés plus précis, cf. théorème \ref{thm:asymp},  corollaire \ref{cor:asym} et proposition \ref{prop:py})
  \begin{enumerate}
  \item Pour tout paramètre $T$,  on a 
$$
\sum_{\of}\int_{G(F)\back G(\AAA)^1}  |k^T_\of(f,g)| \, dg <\infty.
$$
où la somme est prise sur toutes les orbites géométriques.
\item Pour toute orbite $\of$, l'application 
$$T\mapsto \int_{G(F)\back G(\AAA)^1}  k^T_\of(f,g) \, dg $$
est polynomiale en le paramètre $T$.
\item Ce polynôme  est asymptotique à l'intégrale
$$\int_{G(F)\back G(\AAA)^1}    F^G(g,T) \sum_{X\in \of(F)} f(\Ad(g)^{-1}X)\, dg
$$
construite à l'aide de la fonction caractéristique   $F^G(g,T)$ d'un certain compact de $G(F)\back G(\AAA)^1$.
  \end{enumerate}
\end{theoreme}

La série $\sum_{\of}\int_{G(F)\back G(\AAA)^1}  k^T_\of(f,g) \, dg$ est donc absolument convergente et sa somme est \og la formule des traces pour les algèbres de Lie\fg{} (cf. \cite{PH1}). On parle de développement fin car il se fait selon les classes de conjugaison (ici géométriques) alors que dans  \cite{PH1} le développement se fait selon  les classes de conjugaison semi-simples (rationnelles). Pour les corps de fonctions, le groupe $GL(n)$, les orbites nilpotentes et  une fonction test très simple, ce théorème est démontré dans \cite{scfhn}. L'assertion 1 dans le cas des groupes avait été conjecturée par Hoffmann (cf. conjecture 1 \cite{Hoff}). Les méthodes de l'article permettent aussi de démontrer l'énoncé de Hoffmann. L'assertion 3 permet de faire le pont avec les construction d'Arthur dans \cite{ar_unipvar}. On peut définir alors une distribution 
$$J_\of^T(f)=\int_{G(F)\back G(\AAA)^1}  k^T_\of(f,g) \, dg .
$$
Un problème général est de trouver une expression aussi explicite que possible pour cette intégrale en termes de distributions locales plus élémentaires. 
\end{paragr}

\begin{paragr}[Cas de $GL(n)$.] --- Désormais $G=GL(n)$ et $\of$ une orbite \og régulière par blocs\fg{}. Soit $f\in \Sc(\ggo(\AAA)$. On pose 
$$J_\of(f)=J_\of^{T=0}(f).$$
Dans ce cas, il est possible de réécrire l'intégrale qui définit $J_\of(f)$ de sorte qu'on puisse insérer un facteur de convergence. On peut alors permuter l'intégrale et la somme sur $P$ qui définit $k_\of^T(f,g)$. On en déduit une expression pour $J_\of(f)$ en terme d'intégrales orbitales régularisées (cf. théorème \ref{ref:calcul}). En utilisant la théorie des $(G,M)$-familles d'Arthur, on  obtient le théorème suivant (pour les choix implicites de mesures on renvoie au § \ref{S:normalisation}).

\begin{theoreme}(pour une version plus précise, cf. théorème \ref{thm:coef} et  §\ref{S:variante}) Soit $\of$ l'orbite d'un élément de $\ggo$ qui possède un seul bloc de Jordan de taille $r$ avec multiplicité $d$ (on a donc $n=rd$).

Soit $S$ un ensemble fini de places et supposons que $f=f_S\otimes \mathbf{1}^S$ où $\mathbf{1}^S$ est la fonction \og unité\fg{} hors $S$ et $f_S\in\Sc(\ggo(\AAA_S))$. Alors
$$ J_{\of}(f_S\otimes \mathbf{1}^S)= \sum_{(L,\of')}\frac{|W^L|}{|W|}   \tilde{a}^L(S,\of') J_L^G(\of', f_S)
$$
où l'on somme sur les couples $(L,\of')$ formés d'un sous-groupe de Levi semi-standard et d'une $L$-orbite nilpotente $\of'\subset \lgo$ tels que $I_L^G(\of')=\of$. Les groupes $W^L$ et $W$ sont des groupes de Weyl, la distribution $J_L^G(\of', f_S)$ est une intégrale orbitale nilpotente pondérée locale. Enfin le coefficient  $\tilde{a}^L(S,\of')$ admet une expression combinatoire explicite en terme de la fonction $\zeta$ partielle hors $S$ du corps $F$. 
\end{theoreme}

L'expression ci-dessus était formellement connue d'Arthur. L'innovation est que le théorème fournit des formules explicites pour les coefficients  $\tilde{a}^L(S,\of')$. Arthur n'avait donné des expressions que pour l'orbite nulle (c'est le seul cas semi-simple \emph{et} nilpotent).

\end{paragr}

\begin{paragr}[Plan de l'article.] --- La section \ref{sec:notations} rassemble les notations et les notions utiles à la suite. La section \ref{sec:asym}, qui s'applique à tout groupe réductif,  est consacrée à établir un développement fin de la formule des traces. Ensuite l'article se spécialise au cas du groupe $GL(n)$ et d'une orbite régulière par bloc $\of$. La section \ref{sec:IOP} introduit certaines intégrales orbitales régularisées. Celles-ci réapparaissent à la section  \ref{sec:reg} dans  une expression combinatoire pour la contribution $J_\of(f)$. On en déduit, à la section \ref{sec:coef}, le développement de $J_\of(f)$ en termes d'intégrales locales. Dans la section finale \ref{sec:GM}, on démontre  quelques résultats utiles sur les $(G,M)$-familles d'Arthur. 
  
\end{paragr}

\begin{paragr}[Remerciements.] --- Lors de l'élaboration de cet article, j'ai reçu le soutien de l'Institut Universitaire de France et des  projets Ferplay ANR-13-BS01-0012 et  Vargen ANR-13-BS01-0001-01 de l'ANR. Je souhaite également remercier Gérard Laumon : mon intérêt pour la partie unipotente de la formule des traces a été grandement stimulé par notre travail en commun sur le comptage des fibrés de Hitchin (cf. \cite{scfhn}).
\end{paragr}

\section{Notations et rappels}\label{sec:notations}

\begin{paragr} Soit $F$ un corps de nombres. 
\end{paragr}

\begin{paragr} On tâche de suivre les notations d'Arthur. Soit $G$ un groupe algébrique défini sur $F$. On note par la même lettre en minuscule gothique, ici $\ggo$, son algèbre de Lie ; cette convention vaut pour tous les groupes algébriques rencontrés. Soit $N_G$ son radical unipotent.   Soit $X^*(G)$ le groupe des caractères rationnels de $G$ défini sur $F$. L'action adjointe de $G$ sur $\ggo$ est noté $\Ad$. Sauf mention contraire, un sous-groupe de $G$ signifie un sous-groupe algébrique de $G$ défini sur $F$. 
\end{paragr}

\begin{paragr} Supposons $G$ connexe et réductif. Soit $M\subset H$ des sous-groupes de $G$. Soit $\fc^H(M)$ l'ensemble des sous-groupes paraboliques de $G$ inclus dans $H$ et contenant $M$. Supposons de plus que  $M$ est un sous-groupe de Levi de $G$, ce par quoi on entend un facteur de Levi d'un sous-groupe parabolique de $G$. Soit  $\pc^H(M)$ l'ensemble des sous-groupes paraboliques de $G$ inclus dans $H$, dont $M$ est un facteur de Levi.  Tout sous-groupe parabolique $Q\in \fc^H(M)$, possède un unique facteur de Levi noté $M_Q$ qui contient $M$. En particulier, $Q=M_Q N_Q$ est une décomposition de Levi de $Q$. L'image de l'application $Q\in\fc^H(M) \mapsto M_Q$ est notée $\lc^H(M)$. Lorqu'on a  $H=G$, on omet l'exposant $G$ dans la notation.
 \end{paragr}

 \begin{paragr}Soit $M$ un sous-groupe de Levi de $G$. Soit $A_M$ le sous-tore déployé maximal du centre de $G$. Pour tout sous-groupe $H$ de $G$ stable par conjugaison par $A_M$, on note $\Sigma^H_M$ l'ensemble des racines de $A_M$ sur $\hgo$.  Soit  $P\in \pc(M)$ et $Q$ un sous-groupe parabolique contenant $P$.  Soit $\Delta_P^Q$ l'ensemble des racines simples dans $\Sigma^{M_Q\cap N_P}_{M}$. Lorsque $Q=G$, on note simplement $\Delta_P=\Delta_P^G$.   On a $\Delta_P^Q\subset \Delta_P$.
    \end{paragr}
    
\begin{paragr}
      Pour tout $P\in \pc(M)$, soit $a_P=\Hom(X^*(P),\RR)$ et son dual $a_P^*=X^*(P)\otimes_\ZZ\RR$. Soit $\bg \cdot, \cdot \bd$ l'accouplement canonique entre $a_P^*$ et $a_P$. Le morphisme de  restriction $X^*(P)\to X^*(M)$ est bijectif ce qui permet d'identifier $a_P$ à $a_M=\Hom(X^*(M),\RR)$, de même pour les espaces duaux. 

Soit $M\subset L$ des sous-groupes de Levi et $P\in \pc(M)$ et $Q\in \pc(L)$ tels que $P\subset Q$. En utilisant le morphisme de  restriction $X^*(M)\to X^*(A_M)$, on identifie $a_M^*$ à $X^*(A_M)\otimes \RR$. De la sorte on a $\Sigma_M^P\subset a_M^*$. Soit $a_P^{Q,*}$ le sous-espace de $a_P^*$ engendré par $\Delta_P^Q$ ; il ne dépend que de $M$ et $L$, on le note encore $a_M^{L,*}$. On note $\rho_P^Q$ la demi-somme des éléments de  $\Sigma^{M_Q\cap N_P}_{M}$. 
\end{paragr}

\begin{paragr}   
  À chaque $\al\in \Delta_P$ est associée une coracine $\al^\vee\in a_P$ (cf. \cite{ar1}). On note $\Delta_P^{Q,\vee}$ l'ensemble des coracines des éléments de $\Delta_P^Q$. Dualement au morphisme de  restriction $X^*(L)\to X^*(M)$, on a une projection  $a_{M} \to a_{L}$ dont on note $a_M^{L}$ le noyau. L'ensemble $\Delta_P^{Q,\vee}$ forme une base de $a_M^{L}$. Ce dernier est en dualité parfaite avec $a_M^{L,*}$.  L'orthogonal dans $a_M$ de $a_M^{L,*}$ s'identifie via la projection $a_M\to a_L$  à $a_L$. On a donc des décompositions
$$a_M=a_M^L\oplus a_L$$
et de même pour les duaux.
\end{paragr}

\begin{paragr}Soit $\hat{\Delta}_P^Q$ la base des poids de  $a_M^{L,*}$, c'est-à-dire la base duale de $\Delta_P^{Q,\vee}$ : on indexe ses éléments  $(\varpi_\al)_{\al\in \Delta_P^Q}$ avec $\bg \varpi_\be,\al^\vee\bd=\delta_{\al,\be}$ (symbole de Kronecker). Soit  $\hat{\Delta}_P^{Q,\vee}$ la base des copoids, c'est-à-dire la base de  $a_M^{L}$ duale de  $\Delta_P^{Q}$. De même, on indexe ses éléments $(\varpi_\al^\vee)_{\al\in \Delta_P^Q}$. Finalement, pour tout $\Delta\subset a_M$, on note $\ZZ(\Delta)$ le sous-groupe engendré par $\Delta$. 
\end{paragr}

\begin{paragr}
  Lorsqu'on a des groupes notés $P_1\subset P_2$ on note simplement $a_1^2$, $\Delta_1^2$ etc. les objets $a_{P_1}^{P_2}$ et  $\Delta_{P_1}^{P_2}$.
    \end{paragr}

\begin{paragr}[Induction d'orbites.] --- \label{S:ind-LS} Les références pour ce paragraphe sont  \cite{lus-spal} et \cite{Borho}. Soit $\oc=\oc^G$ l'ensemble des orbites adjointes \og géométriques\fg{} de $G$.  Soit $P\subset G$ un sous-groupe parabolique et  $X\in \pgo$. Soit $L=P/N_P$ le plus grand quotient réductif de $P$. Le groupe $L$ agit par adjonction sur son algèbre de Lie $\lgo=\pgo/\ngo_P$. Soit 
$$\pi: \pgo\to \lgo$$
 la projection canonique. 
Soit $\oc^L_X$ l'orbite de $\pi(X)$ sous $L$. 

Soit $I_P^G(X)$, qu'on note encore  $I_P^G(\oc^L_X)$, l'unique élément de $\oc$ tel que l'intersection
\begin{equation}
  \label{eq:intersection}
  I_P^G(X)\cap \pi^{-1}(\oc^L_X)
\end{equation}
soit un ouvert de Zariski  dense dans $\pi^{-1}(\oc^L_X)$. On dit que c'est l'orbite induite. Dans ce cas, l'intersection est exactement la $P$-orbite de $X$.

Identifions le quotient $L$ à un facteur de Levi de $P$. Supposons $X\in \lgo$. Soit $Q\in \pc(L)$. On a alors  $I_Q^G(X)=I_P^G(X)$. On note alors simplement $I_L^G(X)$ cette orbite.

Pour tout sous-groupe parabolique $P\subset Q\subset G$, on peut définir de la même façon une orbite induite $I_P^Q(X)$ qui est une orbite dans $\qgo/\ngo_Q$ (pour l'action adjointe de $Q/N_Q$). L'induction est transitive au sens où l'on a la formule suivante 
\begin{equation}
  \label{eq:transitivite}
  I_Q^G(I_P^Q(X))=I_P^G(X).
\end{equation}

 Le lemme suivant sera utile pour la suite.

  \begin{lemme}\label{lem:induite}
Soit $P$ un sous-groupe parabolique de $G$ de décomposition de Levi $P=M_PN_P$. 
Il existe une famille de polynômes $(\Phi_{i})_{i\in I}$  sur $\mgo_P\times   \ngo_P$ indexée par un ensemble fini $I$  telle que pour toute $M_P$-orbite $\of$ dans $\mgo_P$ il existe $I_\of\subset I$ qui vérifie les deux propriétés suivantes :
\begin{enumerate}
\item pour tout $X\in \of$ et $Y\in   \ngo_P$ on a  $X+Y\in I_P^{P'}(X)$ si et seulement s'il existe $i\in I_\of$ tel que $P_{i}(X,Y)\not=0$. 
\item pour tout $X\in \of$ il existe $i\in I_\of$ tel que le polynôme $P_{\of,i}(X,\cdot)$ soit non identiquement nul sur $\ngo_P$.
\end{enumerate}
\end{lemme}

\begin{preuve}
On fixe une base de $\pgo$. Comme famille de polynômes $(\Phi_i)_{i\in I}$, on prend les déterminants extraits de toute taille de la matrice dans la base choisie de l'endomorphisme de $\pgo$ induit par $\ad(X+Y)$ pour $X\in \mgo_P$ et $Y\in \ngo_P$.

Soit $\of$ une $M_P$-orbite dans $\mgo_P$ et  $X\in \of$. Soit $Y\in  \mgo_{P'}\cap \ngo_P$ et $\oc'$ la $P$-orbite de $X$. Une condition nécessaire et  suffisante pour avoir  $X+Y\in I_P^{P'}(X)$ est que $\oc'$ soit un ouvert dense dans $\oc\oplus \ngo_P$. Pour cela il faut et il suffit qu'on ait égalité entre les espaces tangents en $X+Y$. Ceux-ci étant  $[\pgo,X+Y]$ et $[\mgo_P,X]\oplus\ngo_P$, une condition nécessaire et suffisante pour avoir  $X+Y\in I_P^{P'}(X)$ est donc  l'égalité
\begin{equation}
  \label{eq:CNS}
  [\pgo,X+Y]=[\mgo_P,X]\oplus\ngo_P
\end{equation}
On a toujours $[\pgo,X+Y]\subset[\mgo_P,X]\oplus\ngo_P$. Soit $r_\of=\dim([\mgo_P,X]\oplus\ngo_P)$. L'entier $r_\of$ ne dépend pas du choix de $X\in \of$. La condition nécessaire et suffisante \eqref{eq:CNS} est donc  que le rang de l'endomorphisme  de $\pgo$  induit par $\ad(X+Y)$ soit au moins $r_\of$. L'ensemble $I_\of \subset I$ indexe les déterminants extraits de taille au moins $r_\of$ répond à l'assertion 1. Cet ensemble répond aussi à l'assertion 2 car pour tout $X\in \mgo_P$  il existe au moins un $Y$ tel que $X+Y\in I_P^{P'}(X)$.
\end{preuve}
\end{paragr}

\begin{paragr}
Soit  $\AAA$ l'anneau des adèles de $F$ et $|\cdot|$ la valeur absolue adélique. Pour ne pas alourdir, on ne  la distingue pas dans les notations de la valeur absolue usuelle sur $\CC$.   Soit $G$ un groupe  défini sur $F$.  Soit $G(\AAA)^1$ le noyau de tous les homorphismes $|\chi|$ pour $\chi\in X^*(G)$. On note 
$$[G]=G(F)\back G(\AAA)\text{   et   }   [G]^1=G(F)\back G(\AAA)^1
$$
Les groupes $G(\AAA)$ et $G(\AAA)^1$ est muni d'une mesure de Haar. Les quotients $[G]$ et $[G]^1$ sont munis des mesures quotients par la mesure de comptage.
\end{paragr}

\begin{paragr}
  On suppose désormais  $G$ réductif et connexe.

Soit $P_0$ un sous-groupe parabolique minimal de $G$ muni d'une décomposition de Levi $P_0=M_0N_0$. Soit $P$ un sous-groupe parabolique $P$  standard c'est-à-dire contenant $P_0$. On dispose alors d'une décomposition de Levi standard $P=M_P N_P$ où $M_P$ est l'unique facteur de Levi qui contient $M_0$. 
 
Soit $V$ l'ensemble des places de $F$ et $K=\prod_{v\in V} K_v \subset G(\AAA)$ un sous-groupe compact maximal \og en bonne position\fg{} par rapport à $P_0$. On fixe sur $K$ la mesure de Haar qui donne le volume total $1$. On a la décomposition d'Iwasawa $G(\AAA)=N_P(\AAA) M_P(\AAA)K$. Le logarithme du module des caractères fournit un homomorphisme canonique
$$H_P:P(\AAA) \to a_P.
$$
La décomposition d'Iwasawa induit une application $G(\AAA) \to P(\AAA)/ P(\AAA)\cap K$ qui par composition avec $H_P$ donne une application encore noté $H_P$
$$
H_P:G(\AAA) \to a_P.
$$
\end{paragr}

\begin{paragr}  Soit $A_P^\infty$ la composante neutre du groupe des $\RR$-points du sous-tore $\QQ$-déployé maximal de $\Res_{F/\QQ}A_{M_P}$. Soit $A_{P}^{G,\infty}=A_{P}^\infty\cap G(\AAA)^1$. On a alors une décomposition $G(\AAA)=N_P(\AAA) M_P(\AAA)^1 A_P^\infty K$. On fixe des mesures de Haar sur les différents groupes qui apparaissent de façon à obtenir pour tout $f\in \Cc(G(\AAA)^1)$ la relation 
$$\int_{G(\AAA)^1} f(g), dg =\int_{A_P^{G,\infty}}\int_{M_P(\AAA)^1}\int_{N_P(\AAA)}\int_K \exp(-\bg 2\rho_P,H_P(a)\bd)f(nmak) \, dndmdadk.
$$
\end{paragr}

\begin{paragr}[Théorie de la réduction.] --- Soit $\tau_P=\tau_P^G$ et $\hat{\tau}_P=\hat{\tau}_P^G$ les fonctions caractéristiques respectives des chambres de Weyl ouvertes respectivement aiguë et obtuse dans $a_P^G$. On voit souvent ces fonctions comme des fonctions sur $a_0=a_{P_0}$ par composition avec la projection $a_{0}\to a_P$. Pour alléger les notations, on pose pour $T\in a_0$ et $g\in G(\AAA)$
$$\tau_P(g,T)=\tau_P(H_P(g)-T) 
$$
et
$$\hat{\tau}_P(g,T)=\hat{\tau}_P(H_P(g)-T). 
$$
Si $T=0$, on omet le $T$ dans la notation. On utilise aussi pour $P\subset Q$ les variantes $\tau_P^Q$ et $\hat{\tau}_P^Q$ relatives à $M_Q$ (qu'on voit encore comme des fonctions sur $a_{P_0}$ via la projection  $a_{0}\to a_P^Q$). 

On utilisera les constructions d'Arthur suivantes (cf. \cite{ar1} p.941). Soit $\omega \subset N_0(\AAA)M_0(\AAA)^1$ un ensemble compact et $T_-\in -a_0^+$. Soit $\mathfrak{S}^P(\om,T_-)$ l'ensemble de Siegel formé des 
$$p a k$$
avec $p\in \omega$, $k\in K$ et $a\in A_0^\infty=A_{P_0}^\infty$ tel que $\tau_0^P(a,T_-)=1$. On suppose que $T_-$ et $\omega$ sont tels que  $G(\AAA)=P(F)\mathfrak{S}^P(\om,T_-)$.
Soit  $\mathfrak{S}^P(\om,T_-,T)\subset\mathfrak{S}^P(\om,T_-) $  formé des $x$ tel que pour $\varpi\in \hat{\Delta}_0^{M_P}$
$$\bg \varpi,H_0(x)-T\bd \leq 0.
$$
Soit $F^P(\cdot,T)$ la fonction caractéristique des $g\in G(\AAA)$ pour lesquels il existe $\delta\in P(F)$ tel que $\delta g\in \mathfrak{S}^P(\om,T_-,T)\subset\mathfrak{S}^P(\om,T_-).$ On peut voir cette fonction comme la fonction caractéristique d'un compact de $A_P^\infty N_P(\AAA)M(F)\back G(\AAA)$. 

Soit 
$$
a_0^+
$$
la chambre de Weyl ouverte aiguë positive dans $a_{P_0}^*$. On fixe 
$$T_+\in a_0^+
$$ 
tel que le lemme suivant soit vrai (un tel point existe toujours).

  \begin{lemme}(Arthur, \cite{ar1} lemme 6.4) \label{lem:partition} Soit $P_2$ un sous-groupe parabolique standard. Pour tout $T\in T_++a_0^+$ et tout $g\in P_2(F)\back G(\AAA)$, on a 
$$1=\sum_{P_1\subset P_2} \sum_{\delta\in P_1(F)\back P_2(F) }F^1(\delta g,T) \tau_1^2(\delta g,T).
$$
    \end{lemme}
\end{paragr}

\begin{paragr}
  Soit  $\Sc(\ggo(\AAA))$ l'espace des fonctions complexes de Schwartz-Bruhat sur $\ggo(\AAA)$ munie de sa topologie usuelle (cf. \cite{Weil}). On fixe sur $\ggo$ une forme notée $\bg\cdot,\cdot\bd$, bilinéaire, non dégénérée et invariante par l'action adjointe. On fixe $\psi$ un caractère additif de $F\back \AAA$. On fixe sur $\ggo(\AAA)$ (et de même sur ses sous-algèbre de Lie) la mesure de Haar qui donne le volume $1$ au quotient $\ggo(F)\back\ggo(\AAA)$. On sait alors définir une transformée de Fourier sur  $\Sc(\ggo(\AAA))$.
\end{paragr}

\section{Asymptotique uniforme d'intégrales orbitales tronquées}\label{sec:asym}

\begin{paragr} --- \label{S:kpo}Soit  $f\in \Sc(\ggo(\AAA))$ et $P$ un sous-groupe parabolique  de $G$. On définit pour $g\in N_P(\AAA)M_P(F)\back G(\AAA)$ 
$$
k_P(f,g)=\sum_{X\in \mgo_P(F)}\int_{\ngo_P(\AAA)} f(\Ad(g)^{-1}(X+U)\,dU.
$$

On définit la fonction suivante pour tout $X\in \pgo$ et $\of\in \oc$
$$\xi_{P,\of}(X)=\left\lbrace
  \begin{array}{l}
    1 \text{ si } I_P^G(X)=\of \, ;\\
0 \text{ sinon}.
  \end{array}\right.
$$
Observons que, par définition même, $\xi_{P,\of}(X)$ est invariante par $P$ et ne dépend que de la projection de $X$ sur $\mgo_P$ (selon $\pgo=\mgo_P\oplus\ngo_P$). On a évidemment pour tout $X\in \pgo$
$$\sum_{\of\in \oc} \xi_{P,\of}(X)=1.
$$
On pose pour $g\in N_P(\AAA)M_P(F)\back G(\AAA)$ 
$$
k_{P,\of}(f,g)=\sum_{X\in \mgo_P(F)}\xi_{P,\of}(X) \int_{\ngo_P(\AAA)} f(\Ad(g)^{-1}(X+U)\,dU.
$$
On a donc
$$\sum_{\of\in \oc} k_{P,\of}(f,g)=k_P(f,g).
$$

Soit $T\in a_0$. On définit
\begin{equation}
  \label{eq:kTo}
  k^T_{\of}(f,g)=\sum_{P} \eps_P^G \, \sum_{\delta\in P(F)\back G(F)}  \hat{\tau}_P(\delta g,T)\cdot  k_{P,\of}(f,\delta g)
\end{equation}
où  $\eps_P^G =(-1)^{\dim(a_P^G)}$ et la somme est prise sur l'ensemble fini des sous-groupes paraboliques standard de $G$.

\end{paragr}

\begin{paragr}[Les énoncés.] --- On fixe une norme $\|\cdot\|$ sur l'espace vectoriel $a_0$.

  \begin{theoreme}\label{thm:asymp}
    Soit $\kc$ une partie compacte de $\Sc(\ggo(\AAA))$. Pour tout $\eps>0$, il existe $\eps'>0$ et $c>0$ tel que pour tout $f\in \kc$ et tout $T\in T_++a_0^+$ tel que 
$$
\bg \al,T \bd \geq \eps \| T\|
$$
pour tout $\al\in \Delta$, on a
\begin{equation}
  \label{eq:asymp}
  \sum_{\of\in \oc}\int_{[G]^1}  |  F^G(g,T) \sum_{X\in \of(F)} f(\Ad(g)^{-1}X) - k^T_\of(f,g)| \, dg \leq c\cdot \exp^{-\eps' \|T\|}.
\end{equation}

\end{theoreme}

La preuve du théorème est assez longue et occupe les §§ \ref{S:red1} à \ref{S:preuve3}. Ce théorème a le  corollaire suivant.

\begin{corollaire}\label{cor:asym}
  Pour tout $f\in \Sc(\ggo(\AAA))$ et tout $T\in a_0$ on a 
$$
\sum_{\of\in \oc}\int_{[G]^1}  |k^T_\of(f,g)| \, dg <\infty.
$$
\end{corollaire}
 La preuve du corollaire est donnée au §\ref{S:p-cor}. On peut donc définir une distribution pour tous  $T\in a_0$ et $\of\in \oc$
\begin{equation}
  \label{eq:Jof}
  J_\of^T(f)=\int_{G(F)\back G(\AAA)^1}  k^T_\of(f,g) \, dg.
\end{equation}

On a aussi le résultat suivant.

\begin{proposition}\label{prop:py}
  L'application $T\mapsto J^T_\of(f)$ est polynomiale.
\end{proposition}

La preuve de la proposition est donnée au §\ref{S:py}.

\begin{remarques}
  
  \begin{enumerate}
  \item Dans le cas d'une orbite nilpotente, Arthur a essentiellement montré (du moins dans le cas des groupes) qu'une intégrale  
$$\int_{[G]^1}   F^G(g,T) \sum_{X\in \of(F)} f(\Ad(g)^{-1}X)\, dg
$$ 
est asymptotique à un polynôme en $T$ (cf. \cite{ar_unipvar} théorème 4.2). Seule l'existence du polynôme  est affirmée. Ici on en donne une expression intégrale.
\item Dans le cas des groupes, le corollaire \ref{cor:asym} a été proposé comme conjecture par Hoffmann (cf. \cite{Hoff} conjecture 1). Les méthodes présentées ici doivent aussi fournir une démonstration de la conjecture d'Hoffmann.
  \end{enumerate}
\end{remarques}
\end{paragr}

\begin{paragr}[Preuve de la proposition \ref{prop:py}.] --- \label{S:py}On admet provisoirement le corollaire \ref{cor:asym}. Soit $T_1$ et $T$ deux points de $a_0$. En reprenant les méthodes d'Arthur (cf. par exemple \cite{ar-intro} section 9), on voit qu'il existe pour tout sous-groupe parabolique standard $Q$ un polynôme $p_Q(T_1,\cdot)$ de sorte que 
$$
  J_\of^T(f)=\sum_{P_0\subset Q\subset G} p_Q(T_1,T) \sum_{\of'\in \oc^{M_Q}, I_Q^G(\of')=\of }J^{M_Q,T_1}_{\of'}(f_Q)
$$
où $J^{M_Q,T_1}_{\of'}(f_Q)$ est la distribution qu'on a définie en \eqref{eq:Jof} mais relativement au point $T_1$, au groupe $M_Q$ et à la fonction $f_Q\in\Sc(\mgo_Q(\AAA))$ définie par
$$f_Q(Y)=\int_{\ngo_Q(\AAA)} \int_K f(\Ad(k^{-1})(Y+U))\, dU$$
pour $Y\in \mgo_Q(\AAA)$. La proposition s'en déduit.
\end{paragr}

\begin{paragr}[Preuve du corollaire \ref{cor:asym}.] --- \label{S:p-cor} Comme $F^G( \cdot,T)$ est à support compact sur $[G]^1$, on a évidemment
$$\sum_{\of\in \oc} 
  \int_{[G]^1}    F^G(g,T) \sum_{X\in \of(F)} |f(\Ad(g)^{-1}X)| \, dg <\infty.
$$
Le corollaire \ref{cor:asym} vaut donc pour $T\in T_++a_0^+$ d'après le théorème \ref{thm:asymp}. On fixe $T_1\in  T_++a_0^+$. Soit $T\in a_0$. En reprenant les notations du §\ref{S:py} et  de \cite{trace_inv} section 2 (en particulier la preuve de l'équation (2.3) de \emph{ibid.}), on voit qu'on a la majoration
$$\sum_{\of\in \oc}\int_{[G]^1}  |k^T_\of(f,g)| \, dg \leq \sum_{P_0\subset Q\subset G} \int_{a_Q^G}|\Gamma'_Q(H,T-T_1)|\, dH   \sum_{\of'\in \oc^{M_Q}}  \int_{[M_Q]^1}  |k^{M_Q,T_1}_{\of'}(f_Q,m)| \, dm
$$
où $dH$ est une mesure de Haar sur $a_Q^G$, $\Gamma'_Q(\cdot,T-T_1)$ est une fonction à support compact sur $a_Q^G$ (\cite{trace_inv} lemme 2.1) et $k^{M_Q,T_1}_{\of'}$ est l'expression \eqref{eq:kTo} relative à $M_Q$. Le majorant est fini puisque $T_1\in  T_++a_0^+$.
\end{paragr}

\begin{paragr}[Première réduction dans la démonstration du théorème \ref{thm:asymp}.] --- \label{S:red1}

 Soit $P_2\subset P_3$ des sous-groupes paraboliques standard,  $g\in G(\AAA)$, $\of\in \oc$ et $f\in \Sc(\ggo(\AAA))$. Soit
$$
k_{2,3,\of}(f,g)=\sum_{P_2\subset P\subset P_3} \eps_P^G \cdot k_{P,\of}(f,g)
$$
Arthur a introduit dans \cite{ar1} lemme 6.1  une fonction caractéristique sur $a_0$ notée $\sigma_2^3=\sigma_{P_2}^{P_3}$. Pour $g\in G(\AAA)$, soit
 $$\sigma_2^3(g,T)=\sigma_2^3(H_{P_0}(g)-T).
$$
On rappellera ses propriétés au fur et à mesure de nos besoins. Pour le moment, nous allons montrer que le théorème \ref{thm:asymp} est une conséquence de la proposition suivante.

  \begin{proposition}\label{prop:reduc1}
     Soit $\kc$ une partie compacte de $\Sc(\ggo(\AAA))$. Soit
$$P_1\subset P_2 \subsetneq P_3$$
des sous-groupes paraboliques.
Pour tout $\eps>0$, il existe $\eps'>0$ et $c>0$ tel que pour tout $f\in \kc$ et tout $T\in T_++a_0^+$ tel que 
$$
\bg \al,T \bd \geq \eps \| T\|
$$
pour tout $\al\in \Delta$, on a 
$$\sum_{\of\in \oc}\int_{P_1(F)\back G(\AAA)^1}  |  F^1(g,T_+) F^2(g,T) \tau_1^2(g,T_+) \sigma_2^3(g,T) |k_{2,3,\of}(f,g)| \, dg \leq c\cdot \exp^{-\eps' \|T\|}.
$$
  \end{proposition}
  
Tout d'abord, des manipulations formelles (cf.  pp. 41-43 dans \cite{ar-intro}) montre qu'on a pour tout $\of\in \oc$ et $T\in T_++a_B^+$
$$
k_{\of}^T(f,g)=\sum_{P_2\subset P_3} \sum_{\delta \in P_2(F)\back G(F)} F^2(\delta g,T) \sigma_2^3(\delta g,T) k_{2,3,\of}(f,g).
$$

Pour $P_2=P_3=G$, la contribution est simplement 
$$F^G(g,T) \sum_{X\in \of(F)} f(\Ad(g)^{-1}X).
$$
Pour $P_2=P_3\subsetneq G$ la contribution correspondante est nulle (car alors $\sigma_2^3\equiv 0$ cf. \cite{ar1} lemme 6.1). Pour prouver la majoration \eqref{eq:asymp}, il suffit de la prouver pour toute contribution associée comme ci-dessus à $P_2\subsetneq P_3$. Or, pour $P_2\subsetneq P_3$, on a, par le lemme \ref{lem:partition} (apppliqué au point $T=T_+$)
\begin{eqnarray*}
\int_{[G]^1}   \sum_{\delta \in P_2(F)\back G(F)} F^2(\delta g,T) \sigma_2^3(\delta g,T) |k_{2,3,\of}(f,g)|\, dg\\
=\sum_{P_0\subset P_1\subset P_2} \int_{P_1(F)\back G(\AAA)^1 }F^1(g,T_+)\tau_1^2(g,T_+)   F^2(\delta g,T) \sigma_2^3(\delta g,T) |k_{2,3,\of}(f,g)|\, dg
\end{eqnarray*}
On tombe précisément sur une somme finie d'expressions qui sont majorées par la proposition \ref{prop:reduc1}.
\end{paragr}

\begin{paragr}[Deuxième  réduction dans la démonstration du théorème \ref{thm:asymp}.] ---  Dans ce paragraphe, on va montrer que la proposition \ref{prop:reduc1} résulte  de la proposition suivante.  Pour $P_1\subset P_3$ des sous-groupes paraboliques, soit
  \begin{equation}
    \label{eq:A13}
    A_1^{3,+}=\{a\in A_1^{\infty}\mid \al(a)\geq 1 \ \forall \al\in \Delta_1^3    \}.
  \end{equation}

 \begin{proposition}\label{prop:reduc2}
     Soit $\kc$ une partie compacte de $\Sc(\ggo(\AAA))$,
$$P_1\subset P_2 \subsetneq P_3$$
des sous-groupes paraboliques et 
$$
\al\in \Delta_1^3-\Delta_1^2.
$$
Il existe $c>0$ tel que pour tout $f\in \kc$ et tout $a\in A_1^{3,+}$, 
$$
 \exp(\bg -2\rho_1,H_1(a)\bd) \sum_{\of\in \oc} |k_{2,3,\of}(f,a)|  \leq c\cdot \al(a)^{-1}.
$$
  \end{proposition}

Pour tout $\of\in \oc$, on a, par décomposition d'Iwasawa et nos choix de mesures,
\begin{eqnarray*}
  \int_{P_1(F)\back G(\AAA)^1 }F^1(g,T_+)\tau_1^2(g,T_+)   F^2(\delta g,T) \sigma_2^3(\delta g,T) |k_{2,3,\of}(f,g)|\, dg
\end{eqnarray*}

\begin{equation*}
   = \int_{[N_1]} \int_{[M_1]^1} \int_{A_1^{G,\infty}} \int_K  \exp(\bg -2\rho_1,H_1(a)\bd)F^1(m,T_+) F^2(nma,T)\tau_1^2(a,T_+) \sigma_2^3(a,T) |k_{2,3,\of}(f,nmak)|\, dg
\end{equation*}

En utilisant le fait que $g\mapsto k_{2,3,\of}(f,g)$ et $g\mapsto F^2(g,T)$ sont toutes deux invariantes à gauche par $N_3(\AAA)$, on majore l'expression  ci-dessus par le volume de $[N_3]$ multiplié par
\begin{equation*}
  \int_{[N_1^3]} \int_{[M_1]^1} \int_{A_1^{G\infty}} \int_K  \exp(\bg -2\rho_1,H_1(a)\bd)F^1(m,T_+) F^2(nma,T)\tau_1^2(a,T_+) \sigma_2^3(a,T) |k_{2,3,\of}(f,nmak)|\, dk da dm dn
\end{equation*}

\begin{lemme}\label{lem:cond}
  Soit $a\in A_1^{G,\infty}$, $n\in N_1^3(\AAA)$ et $m\in M_1(\AAA)^1$ et $H\in a_1^G$ tel que $a=\exp(H)$. Sous les trois conditions suivantes

  \begin{itemize}
  \item  $F^2(nma,T)=1$ ;
  \item  $\sigma_2^3(a,T)=1$ ;
  \item $\tau_1^2(a,T_+)=1$,
  \end{itemize}
on a 

\begin{enumerate}
\item Pour tout $\al\in \Delta_1^3-\Delta_1^2$ 
$$\bg\al,H\bd \geq \bg\al,T\bd
$$
\item Pour tout $\al\in \Delta_1^2$ 
$$\bg\al,H\bd \geq \bg\al,T_+\bd
$$
\item soit $H=\sum_{\al \in \Delta_1} H_\al \varpi_\al^\vee$. Pour tout $\al \in \Delta_1^3$, on a  $ 0\leq H_\al$ et si de plus $\al\in \Delta_1^2$ on a 
$$ 0\leq H_\al \leq \frac{\bg \varpi_\al',T\bd }{\bg \varpi_\al',\varpi_\al^\vee\bd}
$$
où  $\varpi_\al'$ est la projection sur $a_1^{2,*}$ de  $\varpi_\al$.
\end{enumerate}
\end{lemme}

\begin{preuve}
  La condition  $F^2(nma,T)=1$ implique que pour tout $\varpi\in \hat{\Delta}_0^2$, on a 
$$\bg \varpi,H_0(nma)-T\bd\leq 0.
$$
Pour  $\varpi\in \hat{\Delta}_1^2\subset  \hat{\Delta}_0^2$, on a
\begin{eqnarray*}
  \bg \varpi,H_0(nma)-T \bd &=&  \bg \varpi,H_0(a)-T \bd\\
&=& \bg \varpi,H-T \bd
\end{eqnarray*}
d'où
\begin{equation}
  \label{eq:positi}
  \bg \varpi,H-T\bd\leq 0.
\end{equation}
La projection de $H-T$ sur $a_1^3$ s'écrit
$$\sum_{\be\in \Delta_2^3} \bg \be,H-T\bd \varpi_\be^\vee+\sum_{\be\in \Delta_1^2} \bg \varpi_\be,H-T\bd \be^\vee.
$$

Soit $\al\in \Delta_1^3-\Delta_1^2$ et $\bar{\al}$ la projection de $\al$ sur $a_2^{3,*}$. On a donc
$$
\bg \al,H-T\bd=  \bg \bar{\al},H-T\bd \varpi_{\bar{\al}}^\vee+\sum_{\be\in \Delta_1^2} \bg \varpi_\be,H-T\bd\bg \al, \be^\vee\bd.
$$
Si $\sigma_2^3(a,T)=1$ le premier terme est positif (cf. lemme 6.1 de \cite{ar1}). Les autres termes le sont aussi en vertu de \eqref{eq:positi} et du fait que la famille $\Delta_1$ est une base \og obtuse\fg{} (c'est bien connu, pour une référence agréable voir \cite{LW} lemme 1.2.6). On obtient ainsi l'assertion 1. L'assertion 2 est une paraphrase de la condition  $\tau_1^2(a,T_+)=1$. On écrit
$$H=\sum_{\al \in \Delta_1} H_\al \varpi_\al^\vee
$$
autrement dit $H_\al=\bg \al,H\bd $ pour $\al\in \Delta_1$. Pour $\al\in \Delta_1^3$, on a  $H_\al \geq 0$ par les assertions 1 et 3 et la positivité de $T$ et $T_+$. Soit $\al \in \Delta_1^2$ et  $\varpi_\al'$  la projection sur $a_1^{2,*}$ de  $\varpi_\al$. On a 
$$ \bg \varpi_\al' ,H\bd=\sum_{\be \in \Delta_1^2} H_\be  \bg \varpi_\al',\varpi_\be^\vee\bd
$$
Comme pour tout $\be,\gamma\in \Delta_1^2$, on a $\bg \gamma,\varpi_\be^\vee\bd \geq 0$, on a aussi $\bg \varpi_\gamma,\varpi_\be^\vee\bd \geq 0$. En particulier  $\bg \varpi_\al',\varpi_\be^\vee\bd\geq 0$. Pour $\al=\beta$, on a 
$$\bg \varpi_\al',\varpi_\al^\vee\bd=\bg \varpi_\al,\varpi_\al^\vee\bd>0.
$$
Par l'inégalité \eqref{eq:positi}, on obtient

$$  H_\al \bg \varpi_\al',\varpi_\al^\vee\bd \leq \bg \varpi_\al',T\bd
$$
ce qui donne la dernière condition.

\end{preuve}

Soit $a_1(T)$ l'ensemble des $H\in a_1^G$ qui vérifie les trois assertions du lemme \ref{lem:cond} et 
$$A_1(T)=\exp(a_1(T)).
$$

Soit  $\Omega_1^3$ et $\Omega$ des sous-ensembles compacts de $N_1^3(\AAA)$ tels que $\Omega_1^3$ s'envoie surjectivement sur $[N_1^3]$ par la projection canonique et pour tout $a\in A_1(T)$ on ait $a^{-1}\Omega_1^3 a\subset \Omega$. De tels compacts existent (cf. la discussion pp. 943-944 de\cite{ar1}). Soit  $\Omega_1\in M_1(\AAA)^1$ un compact tel que la projection canonique $M_1(\AAA)^1\to [M_1]^1$ envoie $\Omega$ surjectivement sur l'ensemble des $m\in   M_1(\AAA)^1$ tel que $F^1(m,T_+)=1$. On est donc ramené à majorer l'expression (à des volumes finis près qu'on peut négliger car indépendants de $T$) 
\begin{equation*}
   \int_{A_1(T)}\exp(\bg -2\rho_1,H_1(a)\bd) \sigma_2^3(a,T)   \sup_{n\in \Omega, m\in \Omega_1, k \in K }|k_{2,3,\of}(f,a nmk)|\, dg
 \end{equation*}
Quitte à remplacer $\kc$ par un compact qui contient les fonctions $f(\Ad(nmk)^{-1}(\cdot))$ pour $n\in \Omega, m\in \Omega_1, k \in K$ et $f\in \kc$, on est réduit à majorer 
\begin{equation*}
   \int_{A_1(T)}\exp(\bg -2\rho_1,H_1(a)\bd)\tau_1^2(a,T_+) \sigma_2^3(a,T)    |k_{2,3,\of}(f,a)|\, da
 \end{equation*}
pour $f\in \kc$. Comme il résulte des définitions, on a  $A_1(T)\subset  A_1^{3,+}$. On est donc en mesure d'utiliser la proposition \ref{prop:reduc2} ou plutôt l'inégalité suivante qui en est un corollaire immédiat : il existe $c>0$ et $e>0$ tels que pour tout $a\in   A_1^{3,+}$
$$
\exp(\bg -2\rho_1,H_1(a)\bd)\sum_{\of\in\oc}   |k_{2,3,\of}(f,a)| \leq c\cdot \prod_{\al\in  \Delta_1^3-\Delta_1^2} \al(a)^{-e}.
$$
La proposition \ref{prop:reduc1} résulte alors du lemme suivant.

\begin{lemme}
  Sous les hypothèses de la proposition \ref{prop:reduc1}, pour tous $e>0$ et tout  $\eps>0$, il existe $\eps'>0$ et $c>0$ tel que pour tout $f\in \kc$ et tout $T\in T_++a_0^+$ tel que 
$$
\bg \al,T \bd \geq \eps \| T\|
$$
pour tout $\al\in \Delta$, on a 
$$
 \int_{A_1(T)} \sigma_2^3(a,T)  \prod_{\al\in  \Delta_1^3-\Delta_1^2} \al(a)^{-e} \, da  \leq c\cdot \exp^{-\eps' \|T\|}.
$$
\end{lemme}

\begin{preuve}
On utilise l'isomorphisme \og exponentielle\fg{} pour identifier la mesure sur $A_1^G$ à une mesure de Haar sur $a_1^G$. Par un changement de variables, il revient au même de majorer
\begin{equation}
  \label{eq:amajorer}
   \int_{a_1(T)} \sigma_2^3(\exp(H),T) \exp(-e\sum_{\al\in  \Delta_1^3-\Delta_1^2} \bg \al,H\bd) \, dH.
 \end{equation}
 Soit $H=H_1+H_2+H_3$ avec
$$H_1=\sum_{\al \in \Delta_1^2} H_\al \varpi_\al^\vee$$
$$H_2=\sum_{\al \in \Delta_1^3-\Delta_1^2} H_\al \varpi_\al^\vee$$
et
$$H_3=\sum_{\al \in \Delta_1-\Delta_1^3} H_\al \varpi_\al^\vee.
$$
Supposons de plus  $H \in a_1(T)$ et $\sigma_2^3(\exp(H),T)=1$. D'après le lemme \ref{lem:cond}, les coefficients de $H_1$ sont majorés par des formes linéaires en $T$. Comme $\sigma_2^3(\exp(H),T)=1$, on a (cf. \cite{ar1} lemme 6.1), on a 
$$\bg \al, H -T \bd \leq 0$$
pour $\al\in \Delta_2-\Delta_2^3$
et
$$\bg \varpi, H -T \bd > 0$$
pour $\varpi\in \hat{\Delta}_3$. En observant que $H_3\in a_3^G$, on voit que $H_3$ est astreint à rester dans un compact de $a_3^G$ qui dépend de $H_1$ et $H_2$. En tenant compte du contrôle des coefficients de $H_1$, on voit qu'il existe un polynôme $P$ en $T$ et $H_\al$ pour $\al\in \Delta_1^3-\Delta_1^2$ de sorte que pour tout $\beta\in  \Delta_1-\Delta_1^3$ on ait
$$|H_\be |\leq P(T,(H_\al)_{\al\in \Delta_1^3-\Delta_1^2}).$$
Ainsi il existe un polynôme $P$ tel que l'expression \eqref{eq:amajorer} soit majorée par 

$$
\int  |P(T,(H_\al)_{\al\in \Delta_1^3-\Delta_1^2})| \prod_{\al \in \Delta_1^3-\Delta_1^2}  \exp(-e H_\al)\,dH
$$
où l'intégrale est prise sur les $H=\sum_{\al \in \Delta_1^3-\Delta_1^2} H_\al \varpi_\al^\vee$  avec $H_\al\geq \bg \al,T\bd$. Il existe $c>0$ et $k>0$ tel que pour $0<e'<e$ on ait 
$$|P(T,(H_\al)_{\al\in \Delta_1^3-\Delta_1^2})| \prod_{\al \in \Delta_1^3-\Delta_1^2}  \exp(-e H_\al)\,dH \leq c \|T\|^k \prod_{\al \in \Delta_1^3-\Delta_1^2}  \exp(-e' H_\al)
$$
Par conséquent \eqref{eq:amajorer} est majorée par
 $$c \|T\|^k \prod_{\al \in \Delta_1^3-\Delta_1^2}  \int^\infty_{\bg \al,T\bd } \exp(-e' H_\al)\, dH_\al= c \|T\|^k \prod_{\al \in \Delta_1^3-\Delta_1^2} (e')^{-1}  \exp(-e' \bg \al,T\bd).
$$
La proposition \ref{prop:reduc1} s'en déduit.
\end{preuve}
\end{paragr}

\begin{paragr}[Troisième  réduction dans la démonstration du théorème \ref{thm:asymp}.] ---  \label{S:reduc3}Dans ce paragraphe, on va montrer que la proposition \ref{prop:reduc2} résulte  de la proposition \ref{prop:reduc3} ci-dessous. On fixe pour toute la suite de la démonstration  du théorème \ref{thm:asymp} une racine 
$$
\al\in \Delta_1^3-\Delta_1^2.
$$
Soit $Q$ le sous-groupe parabolique maximal de $G$ contenant $P_1$ défini par la condition
$$\Delta_1^Q=\Delta_1-\{\al\}.
$$
Comme $\al\notin \Sigma_1^2$, on a aussi $P_2\subset Q$. 

L'application $P\to P\cap M_Q$ induit une application surjective de l'ensemble des sous-groupes paraboliques de $G$ compris entre $P_2$ et $P_3$ sur l'ensemble des sous-groupes paraboliques de  $M_Q$ compris entre $P_2\cap M_Q$ et $P_3\cap M_Q$. Soit $R$ un tel sous-groupe parabolique de $M_Q$. Il a exactement deux antécédents par l'application précédente : l'un, noté $P'$, n'est pas inclus dans $Q$ contrairement à l'autre qui est noté $P$. On a d'ailleurs $P= Q\cap P'$.
Soit
$$
k_{R,\of}(f,g)=\eps_{P'}^G k_{P',\of}(f,g)+\eps_{P}^G k_{P',\of}(f,g).
$$
On a alors

$$
k_{2,3,\of}(f,g)=\sum_{R}k_{R,\of}(f,g)
$$
où la somme porte sur les sous-groupes paraboliques $R$ de $M_Q$  compris entre $P_2\cap M_Q$ et $P_3\cap M_Q$.
La proposition \ref{prop:reduc2} résulte alors de la proposition ci-dessous.

 \begin{proposition}\label{prop:reduc3}      Soit $\kc$ une partie compacte de $\Sc(\ggo(\AAA))$ et $R$ un sous-groupe parabolique de $M_Q$  compris entre $P_2\cap M_Q$  et $P_3\cap M_Q$. Il existe $c>0$ tel que pour tout $f\in \kc$ et tout $a\in A_1^{3,+}$, 
$$
 \exp(\bg -2\rho_1,H_1(a)\bd) \sum_{\of\in \oc} |k_{R,\of}(f,a)|  \leq c\cdot \al(a)^{-1}.
$$
  \end{proposition}
  \end{paragr}

  \begin{paragr}[Démonstration de la proposition  \ref{prop:reduc3}.] ---\label{S:reduc4} On utilise les notations du §\ref{S:reduc3} et les hypothèses de la proposition \ref{prop:reduc3}.

Soit
$$\ngo_P^{P'}=\mgo_{P'}\cap \ngo_P
$$
et
$$
\bar{\ngo}_P^{P'}=\mgo_{P'}\cap \bar{\ngo}_P.
$$
On a 
$$\mgo_{P'}= \bar{\ngo}_P^{P'}\oplus \mgo_P\oplus \ngo_P^{P'}
$$
et
$$\ngo_P=\ngo_P^{P'}\oplus \ngo_{P'}.
$$
Soit 
$$
\tilde{k}_{P',\of}(f,g)=\eps_{P'}^G \sum_{X  \in (\mgo_P\oplus \ngo_P^{P'})(F)}\sum_{Y\in\bar{\ngo}_P^{P'}(F)-\{0\} }\xi_{P',\of}(X+Y) \int_{\ngo_{P'}(\AAA)} f(\Ad(g)^{-1}(X+Y+U)\,dU.
$$

$$
\tilde{k}_{P,\of}(f,g)=-\eps_P^G \sum_{X\in \mgo_P(F)}   \xi_{P,\of}(X) \sum_{Z\in \bar{\ngo}_{P}^{P'}(F)-\{0\}} \int_{\ngo_P(\AAA)} f(\Ad(g)^{-1}(X+U)) \psi(\bg Z,U \bd)\,dU.
$$
et
$$
 \tilde{k}_{R,\of}(f,g)= \sum_{X\in (\mgo_P\oplus\ngo_P^{P'})(F)} (\eps_{P'}^G \cdot \xi_{P',\of}(X)+\eps_P^G \cdot \xi_{P,\of}(X) )  \int_{\ngo_{P'}(\AAA)} f(\Ad(g)^{-1}(X+V))\,dV.
$$

\begin{lemme}\label{lem:dec}
  On a
$$
k_{R,\of}^{}(f,g)=\tilde{k}_{P',\of}(f,g)+\tilde{k}_{P,\of}(f,g)+\tilde{k}_{R,\of}(f,g).
$$
\end{lemme}

\begin{preuve}  Par définition, on a 
$$
k_{R,\of}^{}(f,g)=\eps_{P'}^G \cdot k_{P',\of}(f,g)+\eps_{P}^G \cdot k_{P,\of}(f,g).
$$

On a 
$$\eps_{P'}^G\cdot  k_{P',\of}(f,g)= \tilde{k}_{P',\of}(f,g)+\sum_{X\in (\mgo_P\oplus\ngo_P^{P'})(F)} \eps_{P'}^G \cdot \xi_{P',\of}(X) \int_{\ngo_{P'}(\AAA)} f(\Ad(g)^{-1}(X+V))\,dV.
$$

Il s'agit donc  de voir qu'on a 
$$\eps_{P}^G \cdot k_{P,\of}(f,g)= \tilde{k}_{P,\of}(f,g)+\sum_{X\in (\mgo_P\oplus\ngo_P^{P'})(F)} \eps_{P}^G \cdot \xi_{P,\of}(X) \int_{\ngo_{P'}(\AAA)} f(\Ad(g)^{-1}(X+V))\,dV.
$$

Cette égalité résulte immédiatement de la formule sommatoire de Poisson ci-dessous, relative à la somme rationnelle sur $\ngo_{P}^{P'}(F)$ et valable pour tout $X\in \mgo_P(F)$
\begin{eqnarray}
  \label{eq:Pois}
  \sum_{Y\in \ngo_{P}^{P'}(F)}    \int_{\ngo_{P'}(\AAA)}f(\Ad(g)^{-1}(X+Y+V))\, dV=\\ \nonumber \sum_{Z\in \bar{\ngo}_{P}^{P'}(F)} \int_{\ngo_P(\AAA) } f(\Ad(g)^{-1}(X+U))\psi(\bg Z,U \bd)\, dU.
\end{eqnarray}
\end{preuve}

Le lemme \ref{lem:dec} implique que la  proposition \ref{prop:reduc3} résulte des trois lemmes ci-dessous dont les preuves se trouvent respectivement aux paragraphes \ref{S:preuve1}, \ref{S:preuve2} et \ref{S:preuve3}.

\begin{lemme} \label{lem:maj1}  Pour tout $k \in \NN^*$, il existe $c >0$ tel que pour tout $f\in \kc$ et  tout $a\in A_1^{3,+}$  on ait
 \begin{equation*}
   \exp(\bg -2\rho_1,H_1(a)\bd)  \sum_{\of \in \oc} |\tilde{k}_{P',\of}(f,a)| \leq c \cdot \al(a)^{-k}.
 \end{equation*}
   \end{lemme}

\begin{lemme} \label{lem:maj2}  Pour tout $k \in \NN^*$, il existe $c >0$ tel que pour tout $f\in \kc$ et  tout $a\in A_1^{3,+}$  on ait
 \begin{equation*}
   \exp(\bg -2\rho_1,H_1(a)\bd)  \sum_{\of \in \oc} |\tilde{k}_{P,\of}(f,a)| \leq c \cdot \al(a)^{-k}.
 \end{equation*}
   \end{lemme}

\begin{lemme}\label{lem:reduc3}
  Il existe $c>0$ tel que pour tout $f\in \kc$ et tout $a\in A_1^{3,+}$
\begin{eqnarray*}
 \exp(\bg -2\rho_1,H_1(a)\bd)  \sum_{\of \in \oc} |\tilde{k}_{R,\of}(f,a)| \leq c \cdot \al(a)^{-1}
\end{eqnarray*}
\end{lemme}

\end{paragr}

\begin{paragr} On va utiliser dans la suite le lemme élémentaire suivant.

  \begin{lemme}\label{lem:majoration-base}
L'action par homothétie de $\Gm$ sur $\mathbb{G}_a^n$ induit une action de $\AAA^\times$ sur $\AAA^n$ donc une action de 
$$\RR^\times\hookrightarrow (F\otimes_\QQ \RR)^\times \hookrightarrow \AAA^\times$$
sur $\AAA^n$.

Il existe une constante  $N$ et $c>0$ tel que pour tout

Soit  $f\in \Sc(\AAA^n)$. Pour tout $k\in \NN$, il existe $c>0$ tel que pour tout $a>1$, on ait les assertions suivantes
\begin{enumerate}
\item 
$$\sum_{X\in F^n} f(a X) \leq c \ ;
$$
\item 
$$\sum_{X\in F^n-\{0\}} f(a X) \leq c\cdot a^{-k} \ ;
$$
\item 
$$\sum_{X\in F^n} a^{-n} f(a^{-1} X) \leq c.
$$
\end{enumerate}
\end{lemme}

\begin{preuve}
  Les deux premières assertions sont évidentes. Quant à la troisième, on a par la formule de Poisson
$$
\sum_{X\in F^n} a^{-n} f(a^{-1} X)= \sum_{X\in F^n}  \hat{f}(a X)$$
et l'on est ramené à la première assertion.
\end{preuve}

\end{paragr}

\begin{paragr}[Preuve du lemme \ref{lem:maj1}.] --- \label{S:preuve1} On décompose $\ggo$ en espace propre sous l'action de $A_1$
$$\ggo=\mgo_1\oplus (\oplus_{\beta \in \Sigma_1^G} \ggo_{\beta}).
$$
Pour $\beta \in \Sigma_1^G$, on pose 
$$
d_\beta=\dim_F(\ggo_\be).
$$
On a $d_\be\geq 1$. Tout $X\in \ggo$ s'écrit $X=X_1+\sum_{\al\in \Sigma_1^g} X_\beta$ conformément à cette décomposition. On déduit de \cite{Weil} lemme 5 que pour tout $f$ dans $\kc$, il existe des fonctions positives $\phi_1\in \Sc(\ggo_1(\AAA))$ et $\phi_{\beta}\in \Sc(\ggo_{\beta}(\AAA))$ pour tout $\beta \in \Sigma_1^G$ de sorte que pour tout $X\in \ggo(\AAA)$ on ait
$$
|f(X)|\leq \phi_1(X_1) \prod_{ \be\in \Sigma_1^G  }\phi_{\beta}(X_{\beta}).
$$
L'expression   
$$\exp(\bg -2\rho_1,H_1(a)\bd)  \sum_{\of \in \oc} |\tilde{k}_{P',\of}(f,g) |
$$
 est donc majorée par le produit des expressions suivantes

\begin{equation}
  \label{eq:NP}
  \prod_{\beta\in \Sigma_1^{N_{P'}}} \int_{\ggo_\beta(\AAA)}  \phi_\beta(U_\beta) \,dU_\beta
\end{equation}

\begin{equation}
  \label{eq:1P}
  \prod_{\beta\in \Sigma_1^{ M_{P'}} \cap \Sigma_1^{N_1}} \beta(a)^{-d_\be} \sum_{X_\beta \in \ggo_\beta(F)}\phi_\beta(\beta(a)^{-1}X_\beta) 
\end{equation}

\begin{equation}
  \label{eq:M_1}
  \sum_{X_1 \in \mgo_1(F)}\phi_1(X_1)
\end{equation}

\begin{equation}
  \label{eq:2P}
  \prod_{\beta\in \Sigma_1^{ M_{P}} \cap \Sigma_1^{\bar{N}_1}} \sum_{X_\beta \in \ggo_\beta(F)}\phi_\beta(\beta(a)^{-1}X_\beta) 
\end{equation}

\begin{equation}
  \label{eq:MP-1}
 \sum_{Y \in \bar{\ngo}_P^{P'}(F)-\{0\} } \prod_{\beta\in \Sigma_1^{N_P}\cap \Sigma_1^{M_{P'}}}  \phi_{-\beta}(\beta(a) X_{-\beta}).
\end{equation}
Les expressions \eqref{eq:NP} et \eqref{eq:M_1} sont évidemment finies et indépendantes de $a$. L'expression \eqref{eq:1P} est majorée indépendamment de $a$ (cf. lemme \ref{lem:majoration-base} assertion 3) de même que l'expression \eqref{eq:2P} (cf. lemme \ref{lem:majoration-base} assertion 1). L'expression \eqref{eq:MP-1} est majorée par la somme sur $\gamma \in  \Sigma_1^{N_P}\cap \Sigma_1^{M_{P'}}$ du produit de 

\begin{equation}
  \label{eq:incomplete}
\sum_{X \in \ggo_{-\gamma}(F)-\{0\}} \phi_{-\gamma}(\gamma(a) X)  
\end{equation}
et 
\begin{equation}
  \label{eq:complete}
  \prod_{\beta \in \Sigma_1^{N_P}\cap \Sigma_1^{M_{P'}}-\{\gamma\}} \sum_{X_{\beta}\in \ggo_{-\beta}(F)} \phi_{-\beta}(\beta(a) X_{\beta})
\end{equation}

L'expression \eqref{eq:complete} est majorée indépendamment de $a$ (encore une fois par le  lemme \ref{lem:majoration-base} assertion 1).  En utilisant l'assertion 2 du lemme \ref{lem:majoration-base}, on majore l'expression \eqref{eq:incomplete} par $\gamma(a)^{-k}$ pour tout entier $k\geq 1$.  Comme  $\gamma \in \Sigma_1^{M_{P'}}$, il existe des entiers $n_\beta \geq 0$ de sorte que $\gamma=\sum_{\beta\in \Delta_1^{P'}} n_\be \beta$. Comme, de plus, $\gamma\in \Sigma^{N_P}$, on a $n_\al \geq 1$. En utilisant le fait que $a\in A_1^{3,+}$  et que $\Delta_1^{P'}\subset \Delta_1^3$, il vient
 $$
\gamma(a)^{-k} \leq \al(a)^{-k}
$$
et cela conclut.
\end{paragr}

\begin{paragr}[Preuve du lemme \ref{lem:maj2}.] ---\label{S:preuve2}
 Comme dans le paragraphe \ref{S:preuve1} ci-dessus, il existe des fonctions positives $\phi \in \Sc(\mgo_P(\AAA))$ et $\phi'\in \Sc(\bar{\ngo}_{P}^{P'}(\AAA))$ et $\phi''\in \Sc(\ngo_{P'}(\AAA))$  telles que pour tout $f\in \kc$, tout $X\in \mgo_P(\AAA)$, $V\in \ngo_{P'}(\AAA)$  et tout $Z\in \bar{\ngo}_{P}^{P'}(\AAA)$ on ait
$$
 |\int_{\ngo_{P}^{P'}(\AAA) } f(X+U+V)\psi(\bg Z,U \bd)\, dU| \leq \phi(X)\cdot \phi'(Z)\cdot \phi''(V).
$$

On en déduit (la première majoration utilise un changement de variable évident)

\begin{eqnarray*}
 \exp(\bg -2\rho_1,H_1(a)\bd)  \sum_{\of \in \oc} |\tilde{k}_{P,\of}(f,a)|  \\
\leq   \sum_{X\in \mgo_P(F)} \sum_{Z\in \bar{\ngo}_{P}^{P'}(F)-\{0\}} |(\prod_{\be \in \Sigma_1^{M_P\cap N_1}} \be(a)^{-d_\be})  \int_{\ngo_P(\AAA)} f(\Ad(a)^{-1}(X)+U)) \psi(\bg \Ad(a)^{-1}Z,U \bd)\,dU.\\
\leq [      \sum_{X\in \mgo_P(F)} (\prod_{\be \in \Sigma_1^{M_P\cap N_1}} \be(a)^{-d_\be})  \phi( \Ad(a)^{-1}X)  ]\cdot [\sum_{Z\in \bar{\ngo}_{P}^{P'}(F)-\{0\}} \phi'(\Ad(a)^{-1}Z)] \cdot \int_{\ngo_P^Q(\AAA)}  \phi''(V)\,dV
\end{eqnarray*}
 
En procédant comme dans le paragraphe précédent, on voit que le premier facteur est borné indépendamment de $a$ et que le second peut se majorer par l'expression \eqref{eq:MP-1} pour des fonctions positives $\phi_\be$ bien choisies. Comme précédemment, on arrive à la majoration cherchée. 
\end{paragr}

  \begin{paragr}[Preuve du lemme \ref{lem:reduc3}.] --- \label{S:preuve3}On reprend les notations du §\ref{S:reduc3}. Comme précédemment, il existe des fonctions positives $\phi_1\in \Sc(\mgo_P(\AAA))$, $\phi_2\in \Sc(\ngo_P^{P'}(\AAA))$, $\phi_3\in \Sc(\ngo_{P'}(\AAA))$   de sorte que pour tout $f\in \kc$ et pour tous $X\in \mgo_P(\AAA)$, $Y\in (\ngo_{P}^{P'}(\AAA)$ et $U\in\ngo_{P'}(\AAA) $ on ait
$$|f(X+Y+U)|\leq \phi_1(X) \phi_2(Y)  \phi_3(U)
$$
Il s'ensuit que 
$$\exp(\bg -2\rho_1,H_1(a)\bd) |k_{R,\of}^{(3)}(f,a) |
$$
se majore par le produit de 
$$
 \int_{\ngo_{P'}(\AAA)} \phi_3(U)\, dU
$$  
par la double somme
\begin{equation}
  \label{eq:p1}
  (\prod_{\be \in \Sigma_1^{M_{P}\cap N_1}} \be(a)^{-d_\be}) \sum_{X\in \mgo_P(F)}\phi_1(\Ad(a)^{-1}X)\times
\end{equation}
\begin{equation}
  \label{eq:p2}
  \sum_{Y\in \ngo_P^{P'}(F)} |(\eps_{P'}^G \cdot \xi_{P',\of}(X+Y)+\eps_P^G \cdot \xi_{P,\of}(X) )|   (\prod_{\be \in \Sigma_1^{M_{P'}\cap N_P}} \be(a)^{-d_\be}) \phi_2(\Ad(a)^{-1}Y)   
\end{equation}

\begin{lemme}
  \label{lem:xi}
Soit $X\in \mgo_P$ et $Y\in \ngo_P^{P'}$. Si $X+Y\in I_P^{P'}(X)$ alors  pour toute orbite $\of$
$$
|\eps_{P'}^G \cdot \xi_{P',\of}(X+Y)+\eps_P^G \cdot \xi_{P,\of}(X) |=0.
$$

\end{lemme}
\begin{preuve} Observons que si $X+Y\in I_P^{P'}(X)$ alors on a,  par transitivité de l'induction d'orbites,
$$
I_{P'}^G(X+Y)=I_{P'}^G(I_P^{P'}(X))=I_P^G(X).
$$
Il s'ensuit qu'on a 
$$\xi_{P,\of}(X)=\xi_{P',\of}(X+Y).
$$
Cela conclut car on a toujours
$$
|\eps_{P'}^G \cdot \xi_{P',\of}(X+Y)+\eps_P^G \cdot \xi_{P,\of}(X) |= | \xi_{P',\of}(X+Y)-\xi_{P,\of}(X)  | .
$$
\end{preuve}

On peut donc majorer \eqref{eq:p2} par
\begin{equation}
  \label{eq:p3}
  \sum_{\{Y\in \ngo_P^{P'}(F) \mid X+Y\notin I_P^{P'}(X)\}} (\xi_{P',\of}(X+Y)+\xi_{P,\of}(X) )   (\prod_{\be \in \Sigma_1^{M_{P'}\cap N_P}} \be(a)^{-d_\be}) \phi_2(\Ad(a)^{-1}Y).
\end{equation}
Sommant l'expression \eqref{eq:p3} sur $\of \in \oc$, on obtient 

\begin{equation}
  \label{eq:p4}
  2 \sum_{\{Y\in \ngo_P^{P'}(F) \mid X+Y\notin I_P^{P'}(X)\}}    (\prod_{\be \in \Sigma_1^{M_{P'}\cap N_P}} \be(a)^{-d_\be}) \phi_2(\Ad(a)^{-1}Y).
\end{equation}

 D'après le lemme \ref{lem:induite}, il existe une famille finie de polynômes $(\Phi_i)_{i\in I}$ non tous nuls qui dépendent  de $X\in \mgo_P$ mais dont le degré total de chaque polynôme est borné indépendamment de $X$ telle que pour tout $Y\in \ngo_P^{P'}$ on a
$$X+Y\notin I_P^{P'}(X)$$
si et seulement si pour tout $i\in I$ on a $\Phi_i(Y)=0$. En raisonnant comme dans la preuve de la proposition 5.3.1 de \cite{scfhn}, on démontre qu'il existe une constante $c>0$ telle que pour tout $X\in  \mgo_P(F)$ et tout $a\in A_1^{3,+}$ on ait 
$$ \sum_{\{Y\in \ngo_P^{P'}(F) \mid X+Y\notin I_P^{P'}(X)\}}    (\prod_{\be \in \Sigma_1^{M_{P'}\cap N_P}} \be(a)^{-d_\be}) \phi_2(\Ad(a)^{-1}Y) \leq c\cdot \al(a)^{-1}.$$
On en déduit qu'il existe une constante $c>0$ telle que pour tout $a\in A_1^{3,+}$ la somme
$$\exp(\bg -2\rho_1,H_1(a)\bd) \sum_{\of\in\oc}|k_{R,\of}^{(3)}(f,a) |
$$
soit majorée par le produit de $c \al(a)^{-1}$ par
$$  (\prod_{\be \in \Sigma_1^{M_{P}\cap N_1}} \be(a)^{-d_\be}) \sum_{X\in \mgo_P(F)}\phi_1(\Ad(a)^{-1}X).
$$
Comme précédemment, on montre que cette dernière expression est bornée indépendamment de $a\in  A_1^{3,+}$. Le lemme en résulte.
  \end{paragr}

\section{Intégrales orbitales nilpotentes régularisées}\label{sec:IOP}

\begin{paragr}[Notations] ---\label{S:not-enplus}
Soit   $n\in \NN^*$ et $G=GL(n)$. Soit $T$ le sous-tore maximal de $G$ formé des matrices diagonales. Soit $B\subset G$ le sous-groupe de Borel des matrices triangulaires supérieures et $N_B$ son radical unipotent. Les notations sont celles de la section \ref{sec:notations} à ceci près : les groupes $T$ et $B$ joueront le rôle des groupes $M_0$ et $P_0$. On réserve ces deux dernières lettres pour d'autres objets.
 Soit $P$ un sous-groupe parabolique contenant  $B$. On note  $a_{P,\CC}^*=a_P^*\otimes_\RR\CC$. On définit de même  $a_{P,\CC}^{G,*}$ etc. Soit $a_{P,\CC}^{G,+}$ l'ouvert des $\la\in a_{P,\CC}^{G,*}$ qui vérifient $\Re(\bg \la,\al^\vee\bd)>0$ pour tout $\al\in \Delta_P$. Soit $B\subset P_0\subset P$ des sous-groupes paraboliques. On a une décomposition $a_{0,\CC}^{G,*}=a_{0,\CC}^{P,*}\oplus a_{P,\CC}^{G,*}$ de sorte que tout $\la\in a_{0,\CC}^{G,*}$ s'écrit $\la^P+\la_P$ selon cette décomposition.

\end{paragr}

\begin{paragr}[Normalisation de mesures.] ---\label{S:normalisation}  Pour toute place $v$ de $F$, on munit le corps complété $F_v$ de la mesure de Haar déterminée de la manière suivante :
  \begin{itemize}
  \item si $v$ est réelle, $F_v=\RR$ est muni de la mesure de Lebesgue usuelle ;
  \item si $v$ est archimédienne non réelle , $F_v=\CC$ de la mesure $2dxdy$ ;
  \item si $v$ est non-archmédienne, la mesure est celle qui donne le volume $N(\dc_v)^{-1/2}$ à l'anneau des entiers $\oc_v$ de $F_v$ (où $N(\dc_v)$ est la norme de la différente).
  \end{itemize}

Soit $u_{i,j}$ les coordonnées \og standard\fg{} sur $N_B$. On en déduit la forme volume invariante $\wedge u_{i,j}$  ce qui, avec les choix précédents, permet de  munir $N_B(\AAA)$ et chaque $N_B(F_v)$ d'une mesure de Haar. On a alors $\vol(N_B(\oc_v))=1$ pour presque tout $v$ et
$$\vol(N_B(\AAA)/N_B(F))= 1.$$
Par le même procédé, on munit $\ggo(\AAA)$ et chaque $\ggo(F_v)$ d'une mesure de Haar. On a aussi $\vol(\ggo(\oc_v))=1$ pour presque tout $v$ et $\vol(\ggo(F)\back \ggo(\AAA))= 1.$

Soit $(t_i)_{1\leq i\leq n}$ les   coordonnées \og standard\fg{} sur $T$.  On en déduit la forme volume invariante $\frac{dt_1}{t_1}\wedge \ldots \wedge \frac{dt_n}{t_n}$. Cela permet de munir $T(F_v)$ d'une mesure de Haar pour toute place archimédienne de $F$.  En une place non-archimédienne $v$, on normalise la mesure de Haar sur $T(F_v)$ par 
$$\vol(T(\oc_v))=N(\dc_v)^{-n/2}.$$
Cela munit $T(\AAA_F)$ d'une mesure de Haar. Le morphisme $H_B : T(\AAA)\to a_T$ induit une une suite exacte
\begin{equation}
  \label{eq:suite-ex}
  1 \to T(\AAA)^1  \to  T(\AAA) \to a_T \to   1.
\end{equation}
Au moyen des coordonnées usuelles sur $GL(n)$, on identifie $a_T$ à $\RR^n$ sur lequel on dispose du produit scalaire usuel. L'espace $a_T$ et tous ses sous-espaces sont munis  de la mesure de Haar euclidienne.  On munit $T(\AAA)^1$ de la mesure qui induit sur le quotient  $T(\AAA)^1\back T(\AAA)$ la mesure choisie sur $a_T$.

Soit $K=\prod K_{v}$ le sous-groupe compact maximal de $GL(n,\AAA_F)$ défini par $K_{n,v}=O(n,\RR)$ pour une pace réelle,   $K_{n,v}=U(n,\CC)$ pour une place complexe et  $K_{n,v}=GL(n,\oc_v)$ pour une place non-archimédienne. On met sur $K$ la mesure de Haar qui est le produit des mesures de Haar qui donnent le volume $1$ à chaque $K_v$.

On munit $G(\AAA)$ de la mesure de Haar qui satisfait pour toute $f\in \Cc(G(\AAA))$
$$\int_{G(\AAA)}f(g) \, dg = \int_{T(\AAA)} \int_{N_n(\AAA)} \int_K f(tnk) \, dtdndk.$$
 
Le morphisme $H_G$ induit une suite exacte
 $$
1 \to G(\AAA)^1  \to  G(\AAA) \to a_G \to   1.
$$
et on met sur $G(\AAA)^1$ la mesure qui redonne la mesure euclidienne sur $a_G$.
Cette suite exacte se restreint en une suite exacte
 $$
1 \to T(\AAA)\cap G(\AAA)^1  \to  T(\AAA) \to a_G \to   1
$$
et on munit $T(\AAA)\cap G(\AAA)^1$ de la mesure qui redonne encore la même mesure sur le quotient $a_G$. On a alors  pour toute $f\in \Cc(G(\AAA)^1)$
$$\int_{G(\AAA)^1}f(g) \, dg = \int_{T(\AAA)\cap G(\AAA)^1} \int_{N(\AAA)} \int_K f(tnk) \, dtdndk.$$
Les quotients $[G]$ et $[G]^1$ sont munies des mesures quotients des mesures précédentes par la mesure de comptage sur $G(F)$.
\end{paragr}

\begin{paragr}[Fonctions $\xi$.] ---  Soit $N=[F:\QQ]$ le degré de $F$ sur $\QQ$ et $d_F$ le discriminant de $F$.  Soit $\mathbf{1}=\otimes_v \mathbf{1}_v \in \Sc(\AAA)$ la fonction définie place par place ainsi:
   \begin{itemize}
  \item $\mathbf{1}_v(x)=d_F^{1/2N}\exp(-\pi d_F^{1/N} x^2)$ pour $v$ réelle
\item $\mathbf{1}_v(x)=  d_F^{1/N}\exp(-2\pi d_F^{1/N} |x|^2)$ pour $v$ complexe
  \item $\mathbf{1}_v=N(\dc_v)^{-1/2} \mathbf{1}_{\dc^{-1}_v}$ pour $v$ non-archimédienne où $\mathbf{1}_{\dc^{-1}_v}$ est la fonction caractéristique de la différente inverse.
  \end{itemize}
On vérifie que pour toute place $v$
$$\int_{F_v}\mathbf{1}_v(x)\, dx=1.$$
Soit
$$\xi_v(s)=\int_{F_v^\times} \mathbf{1}_v(x) |x|^s \, dx^{\times}.$$
C'est une intégrale convergente pour $\Re(s)>0$, holomorphe en $s$, qui se prolonge analytiquement au plan complexe en une fonction encoré noté $\xi_v$. Soit
$$\xi(s)=\int_{\AAA^\times} \mathbf{1}(x) |x|^s \, dx^{\times}.$$
C'est une intégrale convergente pour $\Re(s)>1$, holomorphe en $s$, qui se prolonge analytiquement au plan complexe en une fonction encoré noté $\xi$. Pour $\Re(s)>1$, on a la formule du produit
$$\xi(s)=\prod_v \xi_v(s).
$$

\end{paragr}

\begin{paragr}[Fonctions $Z_n$.] --- Pour toute place $v$ de $F$, on introduit la fonction $\mathbf{1}_{\ggo,v}=\mathbf{1}_{\mathfrak{gl}(n),v}\in \Sc(\ggo(F_v)$ qui est défini par $\mathbf{1}_{\ggo,v}(X)=\prod_{1\leq i,j\leq n} \mathbf{1}_v(X_{i,j})$ où $X_{i,j}$ sont les coordonnées de $X$ dans la base canonique. Plus généralement, on définit pour un ensemble quelconque de places,  $\mathbf{1}_{\ggo,S}(g)=\otimes_{v\in S} \mathbf{1}_{\ggo,S} \in \Sc(\ggo(\AAA_S))$ et  $\mathbf{1}_{\ggo}^S(g)=\otimes_{v\notin S} \mathbf{1}_{\ggo}^S \in \Sc(\ggo(\AAA^S))$. Lorsque $S$ est vide, on note  simplement $\mathbf{1}_{\ggo}(g)=\mathbf{1}_{\ggo}^S(g)$.
Soit $s\in \CC$ et 
$$Z_n(s)=\int_{GL(n,\AAA)}   1_{\mathfrak{gl}(n)}  (g) |\det(g)|^{s} dg.
$$
On montre que cette intégrale converge absolument pour $\Re(s)>n$ et que sur l'ouvert  défini par $\Re(s)>0$ on a 
$$Z_n(s+n)=\xi(s+1)\xi(s+2)\ldots \xi(s+n).
$$
ce qui donne \emph{ipso facto} le prolongement méromorphe de $Z_n(s)$ à $\CC$. On introduit aussi la fonction
\begin{equation}
  \label{eq:Ztilde}
  \tilde{Z}_n(s)= (s-n)Z_n(s).
\end{equation}
Alors $\tilde{Z}_n(s+n)$ est holomorphe pour $\Re(s)>-1$. Pour $S$ un ensemble de places, on utilisera les variantes locales
$$Z_{n,S}(s)=\int_{GL(n,\AAA_S)}   1_{\mathfrak{gl}(n),S}  (g) |\det(g)|^{s} dg.
$$
et $Z_n^S(s)=\frac{Z_n(s)}{Z_{n,S}(s)}$. Notons que si $S$ est fini, les fonctions $Z_{n,S}(s+n)$ et
\begin{equation}
  \label{eq:ZtildeS}
  \tilde{Z}^S_n(s+n)= s Z_n^S(s+n).
\end{equation}
sont holomorphes pour   $\Re(s)>-1$.

Avec nos choix de mesure, on a la formule simple suivante.
\begin{equation}
  \label{eq:volume}
  \vol([G]^1)= \vol(a_T^G/\ZZ(\Delta^\vee))\cdot \tilde{Z}_n(n),
\end{equation}
où
\begin{equation}
  \label{eq:covolume}
\vol(a_T^G/\ZZ(\Delta^\vee))=\sqrt{n}.
\end{equation}

On utilisera aussi le lemme suivant.

\begin{lemme}\label{lem:surzeta}
Soit $S$ un ensemble fini de places et $f\in \Sc(\mathfrak{gl}(d,\AAA_S)^r)$. Soit
$$Z(s_1,\ldots,s_r,f)=\int_{GL(d,\AAA_S)^r} f(g_1,\ldots,g_r) \prod_{i=1}^r|\det(g_i)|^{d+s} \, dg
$$
\begin{enumerate}
\item L'intégrale converge absolument pour $\Re(s_i)>-1$.
\item Elle se prolonge en une fonction méromorphe sur $\CC^r$.
\item Soit $\mathbf{1}=\mathbf{1}_{\mathfrak{gl}(d),S}^{\otimes r}$ et $I\subset \{1,\ldots,r\}$. Soit $H_I$ l'hyperplan de $\CC^r$ défini par $s_i=0$ pour tout $i\in I$. Alors pour tout $s=(s_1,\ldots,s_r)$ dans un ouvert de $H_I$, on a 
$$Z(s_1,\ldots,s_r,f)= Z(s_1,\ldots,s_r,\mathbf{1}) \frac{Z((s_i)_{i\notin I},f_I) }{Z((s_i)_{i\notin I},\mathbf{1}_I) }
$$
où 
$$f_I((A_i)_{i\notin I})=\int_{ \mathfrak{gl}(d)(\AAA_S)^I } f(g_1,\ldots,g_r) \, dA_I
$$
où $dA_I$ est la mesure additive sur  $\mathfrak{gl}(d)(\AAA_S)^I$.
\end{enumerate}
\end{lemme}

\begin{preuve}
  Pour l'assertion 1, en majorant $f$ par un produit, on voit qu'il suffit de traiter le cas $r=1$ qui est bien connu (cf. \cite{GJ} Chap.I\,  prop.\,1.1). L'assertion 2 se déduit aussi des méthodes de \emph{ibid.} : on utilise la décomposition d'Iwasawa pour se ramener au cas $d=1$, qui se traite comme dans la thèse de Tate. Il reste à traiter le point 3. La seule observation à faire est qu'il existe une constante $c_S>0$ telle que 
$$dg = c_S |\det(g)|^{-d} d^a g$$
où $dg$ est la  mesure \og multiplicative \og  sur  $GL(d,\AAA_S)$ et celle,\fg{}additive\fg{}, sur $\mathfrak{gl}(d,\AAA_S)$. Ainsi on a dans le cas $r=1$
$$Z(s=0,f)= c_S \int_{  \mathfrak{gl}(d)(\AAA_S)} f(A)\,dA.
$$ 
L'assertion 3 s'en déduit immédiatement.

\end{preuve}

\end{paragr}

\begin{paragr}[Sous-groupe parabolique $P_0$.]\label{S:P0} ---  Soit $r$ et $d$ des entiers $\geq 1$. Soit  $n=r d$ et $G=GL(n)$ avec $n\in \NN^*$.  Soit $P_0$ le sous-groupe parabolique standard de $G$ dont le facteur de Levi standard est isomorphe à $GL(d)^r$. 

Soit $P$ un sous-groupe parabolique de $G$ contenant $P_0$. Il existe un entier $k\geq 1$ et des entiers $n_i\geq1 $ pour $1\leq i\leq k$ de sorte que $r=\sum_{i=1}^k r_i$ et $P$ soit formé des matrices triangulaires supérieures par blocs de la forme
\begin{equation*}
 \left( \begin{array}{ccc}
    A_1 & * & * \\ &\ddots & * \\   & & A_k
  \end{array}
\right)
\end{equation*}
avec $A_i\in GL(r_id)$ (les entrées omise sont des zéros). 

Pour $1\leq i\leq k$, soit 
$$X_{r_id}= \left( \begin{array}{cccc}
    0 & I_d & 0 & 0 \\ &0 &\ddots & 0 \\ & & \ddots & I_d \\  &  & & 0
  \end{array}
\right).
$$
où $I_{d}$ désigne la matrice identité de taille $d$. Le centralisateur de $X_{r_id}$ dans $GL_{r_id}$ est le sous-groupe formé des matrices 

\begin{equation}
  \label{eq:centralisateur}
   \left(
\begin{array}{cccc}
  A_1 & A_2 &  \ldots  & A_r\\
 &A_1  & \ddots & \vdots \\
 
 & & A_1& A_2\\ 
 &  &  & A_1\\
\end{array}\right).
\end{equation}
avec $A_1\in GL(d)$ et $A_i\in \mathfrak{gl}(d)$ pour $2\leq i\leq k$.
Soit $P=MN$ la décomposition de Levi standard de $P$ et
$$
X_P=
 \left( \begin{array}{ccc}
    X_{r_1d} &  &  \\  &\ddots &  \\   &  & X_{r_kd}
  \end{array}
\right)\in \mgo
$$
Lorsque $P=G$, on note simplement $X=X_G$. Soit $M_{X}$ le centralisateur de $X_P$ dans $M$. On a  $M_{X_P}\subset P_0$. Soit
$$P_X=M_{X_P}N_P.
$$
C'est un sous-groupe de $P_0$. C'est encore le centralisateur dans $P$ de la classe de $X+\ngo\in \pgo/\ngo$.

Le groupe $GL(d)^k$ se plonge dans $M_X$ (chaque facteur $GL(d)$ se plonge diagonalement dans $GL(r_id)$ pour $1\leq i \leq k$) et forme un facteur de Levi $L_X$ de $M_X$. On munit $L_X(\AAA)$ de la mesure de Haar donnée par la mesure de Haar sur chaque $GL(d,\AAA)$ (normalisée comme au §\ref{S:normalisation}). Le radical unipotent noté $N_X(\AAA)$ est également muni de la mesure de Haar normalisée par $\vol([N_X])=1$. On en déduit une mesure de Haar \og produit\fg{} sur le groupe unimodulaire $M_X(\AAA)$.
On a des suites exactes  induites respectivement par $H_G$ et $H_P$
\begin{equation}
  \label{eq:exa-MX1}
  1\to M_X(\AAA)\cap G(\AAA)^1 \to  M_X(\AAA) \to a_G\to 1.
\end{equation}
et
\begin{equation}
  \label{eq:exa-MX}
  1\to M_X(\AAA)^1 \to  M_X(\AAA)\cap G(\AAA)^1 \to a_P^G\to 1.
\end{equation}
On munit $M_X(\AAA)\cap G(\AAA)^1$ et  $M_X(\AAA)^1$ des mesures qui redonnent sur les quotients la mesure euclidienne fixée. On observera qu'on a $G_X(\AAA)\cap G(\AAA)^1=G_X(\AAA)^1$. On aura besoin du calcul suivant ; on rappelle qu'on note $[M_X]^1=M_X(F)\back M_X(\AAA)^1.$

\begin{lemme}
  \label{lem:calc-vol} On a 
$$
\vol([M_X]^1)= ( d \tilde{Z}_d(d))^{|\Delta_P^G|+1} \vol(a_T^M/\ZZ(\Delta^{M,\vee}))^{-1}
$$

\end{lemme}

\begin{preuve}
Par nos choix de mesures, on a
$$  \vol([M_X]^1)=\vol([L_X]^1)
$$
où le facteur de Levi $L_X(\AAA)^1$ est muni de la mesure de Haar qui donne sur le quotient $L_X(\AAA)^1\back L_X(\AAA)$ la mesure qu'on obtient en identifiant ce dernier à $a_P^G$ via $H_P$. Le groupe $L_X$ est isomorphe à un produit de groupes $GL(d)$ indexés par les \og blocs\fg{} de $M$. En raisonnant bloc par bloc, on est ramené à traiter le cas  $M=G$. Le volume cherché est alors 
$$\vol([GL_d]^1)=\sqrt{d}\cdot \tilde{Z}_d(d).
$$
multiplié par l'inverse du rapport des  normes
$$
\frac{\|({1}, \ldots,1)\|_{\RR^n}}{\|({1}, \ldots,1)\|_{\RR^d}}=\sqrt{\frac{n}{d}}.
$$
d'où
$$
\vol([G_X]^1)=\frac{d}{\sqrt{n}}  \tilde{Z}_d(d)= d  \tilde{Z}_d(d)  \vol(a_T^G/\ZZ(\Delta^{\vee}))^{-1}  
$$
si l'on tient compte de \eqref{eq:covolume}.
\end{preuve}

\end{paragr}

\begin{paragr}[Polynômes $\theta$ et $\hat{\theta}$.] ---\label{S:theta}
  Suivant Arthur (\cite{trace_inv} p.15), on définit des polynômes de la variable $\la\in a_{P_0}$ pour tout $P_0\subset P$

$$\hat{\theta}_0^P(\la)= \vol(a_0^P/\ZZ(\hat{\Delta}_0^{P,\vee}))^{-1} \prod_{\varpi^\vee \in \hat{\Delta}_0^\vee}\bg \la,\varpi^\vee\bd
$$
et
$$\theta_P(\la)= \vol(a_P^G/\ZZ(\Delta_P^{\vee}))^{-1} \prod_{\al \in \Delta_P} \bg \la,\al^\vee\bd.
$$

\end{paragr}

\begin{paragr}[Intégrale orbitales régularisées.] \label{S:cstr}--- Soit $f\in \Sc(\ggo(\AAA)$. Soit $P\supset P_0$ de décomposition de Levi $P=MN$. Pour tout $\la \in a_{P_0,\CC}^{G,+}$, on pose 

$$
J_{P,X}(f,\la)=\int_{M_X(F) N(F)\back G(\AAA)^1} \exp(-\bg\la,H_0(g)\bd) \cdot\hat{\tau}_P(g) \int_{\ngo_P(\AAA)} f(\Ad(g^{-1})(X_P+U)\,dU\, dg
$$
lorsque l'intégrale converge absolument.

\begin{theoreme}\label{thm:prolongement}
    \begin{enumerate}
  \item Pour tout $\la \in a_{P_0,\CC}^{G,+}$, l'intégrale qui définit $J_{P,X}(f,\la)$ converge absolument et définit une fonction holomorphe de $\la$.
  \item Cette fonction se prolonge en une fonction méromorphe sur $ a_{P_0,\CC}^{G,*}$, encore notée $J_{P,X}(f,\la)$.
  \item La fonction
$$\tilde{J}_{G,X}(f,\la)=\hat{\theta}_0^G(\la) J_{G,X}(f,\la)$$
est holomorphe sur l'ouvert défini par $\Re(\bg \la,\varpi^\vee\bd)>-1$ pour tout $\varpi^\vee\in \hat{\Delta}_0^\vee$.
\item On a 
$$\tilde{J}_{G,X}(f,\la^P)= \hat{\theta}_0^P(\la) J_{P,X}(f,\la)  \theta_P(\la).
$$
   \end{enumerate}
\end{theoreme}

Dans le cas de la fonction particulière $f=\mathbf{1}_{\ggo}$, on peut expliciter les intégrales  $J_{P,X}( f,\la)$.

\begin{proposition}\label{prop:cal-unite}
  On a
  \begin{enumerate}
  \item 
$$J_{P,X}(\mathbf{1}_{\ggo},\la)=\vol([M_X]^1) \theta_P(\la)^{-1} \prod_{\al\in \Delta_0^P} Z_d(d+\frac1d \bg \la^P,\varpi_\al^\vee\bd).
$$
\item 
$$\tilde{J}_{G,X}(\mathbf{1}_{\ggo},\la)=d^{|\Delta_0^G|}\vol(a_0^G/\ZZ(\hat{\Delta}_0^{G,\vee}))^{-1} \cdot \vol([G_X]^1) \cdot \prod_{\al\in \Delta_0^G} \tilde{Z}_d(d+\frac1d \bg \la,\varpi_\al^\vee\bd).
$$
\end{enumerate}
\end{proposition}

La preuve du théorème \ref{thm:prolongement} et de la proposition \ref{prop:cal-unite} occupent le reste de cette section.

\end{paragr}

\begin{paragr}[Preuve de l'assertion 1 du théorème \ref{thm:prolongement}] ---\label{S:pr1}
  On commence par faire une série de manipulations formelles qu'on justifiera ensuite. Soit
$$f_P(Y)=\int_K\int_{\ngo_P(\AAA)} f(\Ad(k^{-1})(Y+U))\,dU. 
$$
Par décomposition d'Iwasawa, on a 
$$
J_{P,X}(f,\la)= \int_{M_X(F) \back M(\AAA)\cap G(\AAA)^1}  \exp(-\bg\la,H_0(m)\bd) \cdot\hat{\tau}_P(H_P(m)) f_P(\Ad(m)^{-1}(X_P))\,dm.
$$
On commence par calculer pour $m\in  M(\AAA)\cap G(\AAA)^1$ l'intégrale
$$\int_{M_X(F) \back M_X(\AAA)\cap G(\AAA)^1}  \exp(-\bg\la,H_0(h m)\bd) \cdot\hat{\tau}_P(H_P(h m)) \,dh.
$$
Observons qu'on a $M_X\subset M\cap P_0$ et  donc $H_0(h m)=H_0(h)+H(m)$. En utilisant la suite exacte \eqref{eq:exa-MX}, on voit que l'intégrale s'écrit encore
$$ \exp(-\bg\la,H_0(m)\bd)\vol([M_X]^1)\int_{a_P^G}  \exp(-\bg\la,H)\bd) \cdot\hat{\tau}_P(H+H_P( m)) \,dH.
$$
Un calcul usuel montre que cette intégrale converge pour $\Re(\bg \la,\al^\vee\bd)>0$ pour tout $\al\in \Delta_P$ et qu'elle vaut sur cet ouvert
$$
\exp(-\bg\la^P,H_0(m)\bd)\vol([M_X]^1) \theta_P(\la)^{-1}.
$$ 
On voit donc qu'il suffit de prouver la convergence absolue de l'intégrale
\begin{equation}
  \label{eq:simpl1}
  \int_{M_X(\AAA)\back M(\AAA)}  \exp(-\bg\la^P,H_0(m)\bd) f_P(\Ad(m)^{-1}X_P)\, dm.
\end{equation}
On utilise la décomposition d'Iwasawa de $M(\AAA)$ :
$$M(AAA)= N_0^P(\AAA) M_0(\AAA) K_M $$
avec $N_0^P=N_0\cap M$ et $K_M=K\cap M(\AAA)$. On a une décomposition de Levi pour $M_X(\AAA)$
$$M_X(\AAA)= N_{0,X_P}^P(\AAA)  M_{0,X_P}(\AAA)$$
où l'ajout de l'indice $X_P$ signifie qu'on prend le centralisateur. En observant que la fonction $f_P$ est $K\cap M(\AAA)$-invariante par adjonction et en tenant compte de $\vol(K_M)=1$, on voit que l'intégrale \eqref{eq:simpl1} s'écrit
\begin{equation}
  \label{eq:new-int}
  \int_{   M_{0,X_P}(\AAA)\back  M_0(\AAA)}  \exp(-\bg\la^P+2\rho_0^P,H_0(m)\bd)  \int_{  N_{0,X_P}^P(\AAA)  \back N_0^P(\AAA) }f_P(\Ad(m)^{-1}\Ad(n)^{-1}X_P)\, dn dm.
\end{equation}
Soit $\ngo^{P,\der}_0=[\ngo^P_0,\ngo_0^P]$. L'application $N_{0,X_P}^P\back N_{0}^0 \to \ngo^{P,\der}_0$ donnée par $n\mapsto \Ad(n^{-1})X_P -X_P$ est un isomorphisme qui induit une bijection au niveau des ensembles adéliques qui préservent la mesure (cf. \cite{scuft} lemme 7.3.3). Ainsi on voit que \eqref{eq:new-int} est égale à 
\begin{equation}
  \label{eq:new-int2}
  \int_{   M_{0,X_P}(\AAA)\back  M_0(\AAA)}  \exp(-\bg\la^P+2\rho_0^P,H_0(m)\bd)  \int_{ \ngo^{P,\der}_0(\AAA) }f_P(\Ad(m)^{-1}(X_P+U))\, dU dm.
\end{equation}
 
 On a un supplémentaire naturel de $\ngo^{P,\der}_0$ dans $\ngo_0^P$ qui est stable par l'action adjointe de $M_0$ à savoir la somme des espaces poids 
$$\tilde{\ngo}_0^P=\oplus_{\al\in \Delta_0^P}\ngo_\al.
$$
Chaque espace $\ngo_\al$ s'identifie naturellement à $ \mathfrak{gl}(d)$. Soit 
$$ \iota: \mathfrak{gl}(d)^{\Delta_0^P}\to \tilde{\ngo}_0^P.
$$
l'isomorphisme qui s'en déduit. On a  $X_P\in \ngo^{P,\der}_0$. L'application $m\mapsto \iota^{-1}(\Ad(m^{-1})X_P)$ identifie le quotient $M_{0,X_P}(\AAA)\back  M_0(\AAA)$ à $GL(d,\AAA)^{\Delta_0^P}$ et cette identification est compatible au mesure. Si l'on fait dans \eqref{eq:new-int2} le changement de variable $U\mapsto \Ad(m)^{-1}U$ et cette identification, on obtient que l'intégrale \eqref{eq:new-int2} est égale à 
$$\int_{GL(d,\AAA)^{\Delta_0^P}} \tilde{f}_{P}(\iota((g_\al)_{\al\in \Delta_0^P})  \prod_{\al\in \Delta_0^P} |\det(g_\al)|^{d+\frac1d \bg\la^P,\varpi_\al^\vee\bd  } dg_\al
$$
où l'on a défini la fonction $\tilde{f}_{P}\in \Sc(  \mathfrak{gl}(d)(\AAA)^{\Delta_0^P})$ par
$$\tilde{f}_{P}(Y)=\int_{ \ngo^{P,\der}_0(\AAA)   }  f_P(\iota(Y)+U)\,dU  $$
On peut trouver des fonctions $\phi_\al\in \Sc(\mathfrak{gl}(d)(\AAA))$ telles que 
$$ \tilde{f}_{P}(\iota((g_\al)_{\al\in \Delta_0^P}) ) \leq \prod_{\al \in \Delta_0^P} \phi_\al(g_\al).
$$
La convergence voulue se déduit de celle de chaque intégrale
$$\int_{GL(d,\AAA)}   \phi_\al(g) |\det(g)|^{s+d} dg
$$
pour $\Re(s)>0$. Cette convergence se voit aisément par une décomposition d'Iwasawa qui ramène l'intégrale en une intégrale sur le tore diagonal. On montre également que la convergence est uniforme pour $s$ dans un compact inclus dans $\Re(s)>0$. On en déduit l'assertion 1 du théorème \ref{thm:prolongement}.
\end{paragr}

\begin{paragr}[Preuve de la proposition \ref{prop:cal-unite} et preuve du  \ref{thm:prolongement} pour la fonction $\mathbf{1}_{\ggo}$.] --- Avec les notations du paragraphe précédent, on a  $(\mathbf{1}_{\ggo})_P=\mathbf{1}_{\mgo}$. De même, on a $(\widetilde{  \mathbf{1}_{\ggo}})_P=\otimes_{\al\in \Delta_0^P} 1_{\mathfrak{gl}(d)}  $. En reprenant les calculs précédents, on arrive à 
$$J_{P,X}(\mathbf{1}_{\ggo},\la)=\vol([M_X]^1) \theta_P(\la)^{-1} \prod_{\al\in \Delta_0^P} Z_d(d+\frac1d \bg \la^P,\varpi_\al^\vee\bd).
$$
Cela donne l'assertion 1 de la   proposition \ref{prop:cal-unite} mais aussi l'assertion 2 (le prolongement méromorphe) pour cette fonction particulière. On a donc
$$\tilde{J}_{G,X}(\mathbf{1}_{\ggo},\la)=d^{|\Delta_0^G|}\vol(a_0^G/\ZZ(\hat{\Delta}_0^{G,\vee}))^{-1} \cdot \vol([G_X]^1) \cdot \prod_{\al\in \Delta_0^G} \tilde{Z}_d(d+\frac1d \bg \la,\varpi_\al^\vee\bd)
$$
qui est bien holomorphe pour $\la$ tel que $\Re(\bg \la^P,\varpi_\al^\vee\bd)>-d$ pour tout $\al\in \Delta_0^P$. On a  alors
\begin{eqnarray*}
  \tilde{J}_{G,X}(\mathbf{1}_{\ggo},\la^P)&=& d^{|\Delta_0^G|}\vol(a_0^G/\ZZ(\hat{\Delta}_0^{G,\vee}))^{-1} \cdot \vol([G_X]^1) \cdot\tilde{Z}_d(d)^{|\Delta_P|}  \prod_{\al\in \Delta_0^P} \tilde{Z}_d(d+\frac1d \bg \la^P,\varpi_\al^\vee\bd) \\
&=& d^{|\Delta_P|} \frac{\vol(a_0^P/\ZZ(\hat{\Delta}_0^{P,\vee}))}{\vol(a_0^G/\ZZ(\hat{\Delta}_0^{G,\vee}))}\cdot \frac{ \vol([G_X]^1)}{ \vol([M_X]^1)}\cdot\tilde{Z}_d(d)^{|\Delta_P|}  \cdot\hat{\theta}_0^P(\la) J_{P,X}( \mathbf{1}_{\ggo} ,\la)  \theta_P(\la)
\end{eqnarray*}

\begin{lemme}\label{lem:volGvolM}On a 
$$   \frac{\vol(a_0^P/\ZZ(\hat{\Delta}_0^{P,\vee}))}{\vol(a_0^G/\ZZ(\hat{\Delta}_0^{G,\vee}))}\cdot \frac{ \vol([G_X]^1)}{ \vol([M_X]^1)}\cdot (d \tilde{Z}_d(d))^{|\Delta_P|}=1.
$$
\end{lemme}

\begin{preuve}
Par nos choix de mesures, on a 
$$\vol(a_T^P/\ZZ(\hat{\Delta}_B^{P,\vee})) \cdot \vol(a_T^P/\ZZ({\Delta}_B^{P,\vee}))=1.
$$
En utilisant les suites exactes pour $P$ et $P=G$
$$0\to a_T^{P_0} \to a_T^P \to a_{P_0}^P \to 0
$$
on a 
  \begin{eqnarray*}
    \frac{\vol(a_0^P/\ZZ(\hat{\Delta}_0^{P,\vee}))}{\vol(a_0^G/\ZZ(\hat{\Delta}_0^{G,\vee}))}=\frac{\vol(a_T^P/\ZZ(\hat{\Delta}_B^{P,\vee}))}{\vol(a_T^G/\ZZ(\hat{\Delta}_B^{G,\vee}))}=\frac{\vol(a_T^G/\ZZ({\Delta}_B^{G,\vee}))}{\vol(a_T^P/\ZZ({\Delta}_B^{P,\vee}))}
  \end{eqnarray*}
Le lemme résulte alors du lemme \ref{lem:calc-vol}.
\end{preuve}

Le théorème \ref{thm:prolongement} est donc complètement prouvé pour la fonction particulière $\mathbf{1}_{\ggo}$.
\end{paragr}

\begin{paragr}[Fin de la preuve du théorème \ref{thm:prolongement}.] ---    Soit $f\in \Sc(\ggo(\AAA)$ et $S$ un ensemble fini de places telles que $f=f_S\otimes \mathbf{1}_\ggo^S$ avec $f_S\in \Sc(\ggo(\AAA_S)$.
En reprenant les calculs précédents, on obtient 

$$J_{P,X}(f,\la)= \frac{J_{P,X}^S(f_S,\la)}{J_{P,X}^S(\mathbf{1}_{S,\ggo},\la)} \cdot J_{P,X}(\mathbf{1},\la).
$$
où l'on pose
$$J_{P,X}^S(f_S,\la)=\int_{GL(d,\AAA_S)^{\Delta_0^P}} \widetilde{(f_S)}_{P}(\iota((g_\al)_{\al\in \Delta_0^P})  \prod_{\al\in \Delta_0^P} |\det(g_\al)|^{d+\frac1d \bg\la^P,\varpi_\al^\vee\bd  } dg_\al
$$
où  $\widetilde{(f_S)}_{P}$ est l'analogue local évident de $f\mapsto \tilde{f}_P$. L'assertion 2 du théorème \ref{thm:prolongement} résulte alors   de l'assertion 2 du lemme \ref{lem:surzeta} ci-dessous et  de l'assertion analogue pour la fonction particulière $\mathbf{1}_\ggo$.
On en déduit également
$$\tilde{J}_{G,X}(f,\la)= \frac{J_{G,X}^S(f_S,\la)}{J_{G,X}^S(\mathbf{1}_{S,\ggo},\la)} \cdot \tilde{J}_{G,X}(\mathbf{1},\la),
$$
ce qui donne l'assertion 3 du théorème \ref{thm:prolongement} (cf. assertion 1 du  lemme \ref{lem:surzeta}). En utilisant l'assertion 3 du lemme \ref{lem:surzeta} et les expressions obtenues au §\ref{S:pr1}, on montre que, pour $\la$ dans l'ouvert de convergence, on a 

$$\frac{J_{G,X}^S(f_S,\la^P)}{J_{G,X}^S(\mathbf{1}_{S,\ggo},\la^P)}= \frac{J_{P,X}^S(f_S,\la)}{J_{P,X}^S(\mathbf{1}_{S,\ggo},\la)}.
$$
Il s'ensuit qu'on a aussi
\begin{eqnarray*}
  \tilde{J}_{G,X}(f,\la^P)&=& \frac{J_{G,X}^S(f_S,\la^P)}{J_{G,X}^S(\mathbf{1}_{S,\ggo},\la^P)} \cdot \tilde{J}_{G,X}(\mathbf{1},\la^P)\\
&=&  \frac{J_{G,X}^S(f_S,\la^P)}{J_{G,X}^S(\mathbf{1}_{S,\ggo},\la^P)}\cdot \hat{\theta}_0^P(\la) J_{P,X}(\mathbf{1}_\ggo,\la)  \theta_P(\la)\\
&=& \frac{J_{P,X}^S(f_S,\la)}{J_{P,X}^S(\mathbf{1}_{S,\ggo},\la)} \cdot \hat{\theta}_0^P(\la) J_{P,X}(\mathbf{1}_\ggo,\la)  \theta_P(\la)\\
&=& \hat{\theta}_0^P(\la) J_{P,X}(f,\la)  \theta_P(\la)
\end{eqnarray*}
d'abord sur un ouvert non vide puis partout par prolongement analytique.
\end{paragr}

\section{Contribution des orbites régulières par bloc pour $GL(n)$}\label{sec:reg}

\begin{paragr}
  Les notations sont celles de la section \ref{sec:IOP}. On dispose de l'élément nilpotent  $X$ \og régulier par bloc\fg{} défini en  §\ref{S:P0}.  Soit $\of$ l'orbite de $X$ et $J_\of^T$ la distribution  sur $\ggo(\AAA)$ définie en \eqref{eq:Jof}. On pose
$$J_\of(f)= J_\of^{T=0}(f)$$
pour tout $f\in \Sc(\ggo(\AAA))$. Soit $W$ le groupe de Weyl de $(G,T)$. On l'identifie au sous-groupe des matrices de permutations. Pour $M\in \lc(T)$, on a une version relative $W^M$. Rappelons qu'on dispose du groupe parabolique $P_0=M_0N_0$ (cf. §\ref{S:P0}). Soit $W_0=\Norm_{W}(M_0)/W^{M_0}$. Pour tout $L\in \lc(M_0)$, on a aussi  $W_0^L=\Norm_{W^L}(M_0)/W^{M_0}$. Le groupe $W$ agit naturellement sur $a_T$ et son dual. Le groupe $W_0^L$ agit naturellement sur  $a_0^L$ et son dual.

  \begin{theoreme}
    \label{ref:calcul}
Soit $f\in \Sc(\ggo(\AAA))$. La fonction méromorphe définie par
$$
\frac{1}{|W_0|} \sum_{w\in W_0} \tilde{J}_{G,X}(f,w\la) \theta_{P_0}(w\la)^{-1}
$$
est holomorphe en $\la=0$ et sa valeur en ce point est égale à $J_\of(f)$.
  \end{theoreme}

  \begin{remarque}
    Dans le cas des groupes et de  l'orbite unipotente régulière, une telle formule était connue de Lapid et Finis comme je l'ai appris lors d'un exposé de Finis à Orsay en 2010. 
  \end{remarque}
\end{paragr}

\begin{paragr}[Preuve du théorème \ref{ref:calcul}.] --- Soit $f\in \Sc(\ggo(\AAA))$. On pose pour $g\in P_0(F)\back G(\AAA)$

$$k_{P,X}(f,g) =\left\lbrace
  \begin{array}{l}
  \displaystyle  \sum_{\delta\in P_X(F)\back P_0(F)}  \int_{\ngo_P(\AAA)}  f(\Ad(\delta g)^{-1}(X+U))\,dU  \text{   si  } P_0\subset P\\
0 \text{  sinon}.
\end{array}
\right.
$$
Soit
$$
k_X(f,g)=\sum_{P_0\subset P\subset G}  \eps_P^G\cdot  \hat{\tau}_P( g) \cdot k_{P,X}( g) 
$$
C'est une variante des constructions du §\ref{S:kpo}.  On a 
 $$k_{P,\of}(f,g)=\sum_{\delta\in P_0(F)\back P(F)}k_{P,X}(f,\delta g)
$$
et
\begin{equation}
  \label{eq:kokX}
  k_\of(f,g)=\sum_{\delta\in P_0(F)\back G(F)} k_X(f,\delta g)
\end{equation}
valables pour tout $g\in G(\AAA)$.

 On utilise alors le théorème de convergence suivant.

\begin{theoreme}\label{thm:cv}
Soit $f\in \Sc(\AAA)$. Pour $\la=0$ ou pour $\la\in a_{P_0,\CC}^{G,+}$ l'intégrale suivante converge.

$$
\int_{P_0(F)\back G(\AAA)^1}   | \exp(-\bg \la, H_{P_0}(g)\bd)\cdot   k_{X}(f, g) |\,dg<\infty
$$
\end{theoreme}

\begin{preuve}
  Considérons le cas  $\la\in a_{P_0,\CC}^{G,+}$. Formellement, on a 
\begin{eqnarray*}
  \int_{P_0(F)\back G(\AAA)^1}    \exp(-\bg \la, H_{P_0}(g)\bd)\cdot   k_{X}(f, g) \,dg&=& \sum_{P_0\subset P\subset G} \eps_P^G \cdot J_{P,X}(f,\la).
\end{eqnarray*}
avec les notations du §\ref{S:cstr}. D'après le théorème \ref{thm:prolongement}, l'intégrale qui définit $ J_{P,X}(f,\la)$ converge absolument. La convergence absolue de l'intégrale considérée s'ensuit. Le cas $\la=0$ est traité à part au §\ref{S:fin-de-p}.
\end{preuve}

On en déduit le lemme.

\begin{lemme} \label{lem:lim}
  Soit  $\la\in a_{P_0,\CC}^{G,+}$. On a  
$$\lim_{t\to 0^+} \int_{P_0(F)\back G(\AAA)^1}   \exp(-\bg t\la, H_{P_0}(g)\bd) \cdot   k_{X}(f, g) \,dg= J_\of(f).
$$
\end{lemme}

\begin{preuve}
  Pour tout $t$ dans l'intervalle $[0,1]$, on a pour tout $g\in P_0(F)\back G(\AAA)^1$
$$
|\exp(-\bg t\la, H_{P_0}(g)\bd)| \leq (1-t) +t \exp(-\bg \Re\la), H_{P_0}(g)\bd) \leq 1 +  \exp(-\bg \Re(\la), H_{P_0}(g)\bd)
$$
En particulier,
$$|\exp(-\bg t\la, H_{P_0}(g)\bd)|\leq  |k_{X}(f, g)| (1 +  \exp(-\bg \Re(\la), H_{P_0}(g)\bd))
$$
et d'après le théorème  \ref{thm:cv} la fonction majorante est absolument convergente. Il s'ensuit que l'intégrale considérée est une fonction continue de $t\in [0,1] $. La limite cherchée existe et vaut
$$\int_{P_0(F)\back G(\AAA)^1}     k_{X}(f, g) \,dg= \int_{[G]^1} \sum_{\delta\in P_0(F)\back G(F)}  k_{X}(f, \delta g) \,dg
$$
et cette dernière intégrale est bien $J_\of(f)$ (cf. \eqref{eq:kokX}).
\end{preuve}

La proposition ci-dessous donne une expression pour $J_\of(f)$.  Combinée avec des résultats généraux sur les $(G,M)$-familles d'Arthur (cf. proposition \ref{prop:cp} et corollaire \ref{cor:cp}), elle implique le théorème \ref{ref:calcul}.

\begin{proposition}\label{prop:calcul}
  L'expression 
$$\sum_{P_0\subset P\subset G} \eps_P^G \cdot \hat{\theta}_0^P(\la)^{-1} \tilde{J}_{G,X}(f,\la^P) \theta_P(\la)^{-1}
$$
est holomorphe en $\la=0$ et sa valeur en ce point est égale à $J_{\of}(f)$.
\end{proposition}

\begin{preuve}
La fonction considérée est méromorphe. L'holomorphie en $\la=0$ résulte de celle de  $\tilde{J}_{G,X}(f,\la)$ au voisinage de $0$ (théorème \ref{thm:prolongement} assertion 3) et d'un théorème général (cf. théorème 2.6.4 de \cite{scfhn}). Le seul point restant est le calcul de la valeur en $\la=0$. Compte tenu du lemme   \ref{lem:lim}, il suffit de voir que la fonction considérée coïncide sur l'ouvert $a_{P_0,\CC}^{G,+}$ avec l'intégrale 
$$\int_{P_0(F)\back G(\AAA)^1}    \exp(-\bg \la, H_{P_0}(g)\bd)\cdot   k_{X}(f, g) \,dg.
$$
Or on a vu dans la preuve du théorème \ref{thm:cv} que cette intégrale est égale pour $\la\in a_{P_0,\CC}^{G,+}$ à 
$$\sum_{P_0\subset P\subset G} \eps_P^G \cdot J_{P,X}(f,\la).$$
On utilise alors le théorème \ref{thm:prolongement} assertion 4 pour conclure.
\end{preuve}

\end{paragr}

\begin{paragr}[Fin de la preuve du théorème \ref{thm:cv}.] ---\label{S:fin-de-p} Rappelons qu'il reste à prouver la convergence de 
$$
\int_{P_0(F)\back G(\AAA)^1}   | k_{X}(f, g) |\,dg.
$$
Dans le cas qu'on considère $G=GL(n)$,  le lemme \ref{lem:partition} vaut pour le point $T=0$. On note alors $F(g)=F(g,T=0)$, $\tau_P(g)=\tau(g,T=0)$ etc. Pour alléger encore, on omet aussi $f$ dans les notations. On a alors  pour $g\in G(\AAA)$
  $$k_X(g)=\sum_{B\subset P\subset G} \eps_P^G\cdot  \hat{\tau}_P( g) \cdot k_{P,X}( g)
$$
$$
=\sum_{B\subset P\subset G} \eps_P^G\big[\sum_{B \subset P_1 \subset P} \sum_{\nu\in P_1(F)\back P(F)} F^1(\nu g) \tau_1^P(\nu g) \big]  \hat{\tau}_P(g) \cdot k_{P,X}(g)
$$
$$
=\sum_{P_1\subset P_2} \sum_{\nu\in P_1(F)\back G(F)}F^1(\nu g) \sigma_1^2(\nu g) \sum_{\{P \mid \nu \in P(F) \text{ et }P_1\subset P\subset P_2 \}}  \eps_P^G\cdot k_{P,X}(g).
$$
Par le changement de variable $\nu'=\nu\delta$, on a
$$ \sum_{\delta \in P_0(F)\back G(F)}   |\sum_{B\subset P\subset G} \eps_P^G\cdot  \hat{\tau}_P(\delta g) \cdot k_{P,X}(\delta g) |$$
$$= \sum_{\delta \in P_0(F)\back G(F)}  |\sum_{P_1\subset P_2} \sum_{\nu\in P_1(F)\back G(F)}F^1(\nu g) \sigma_1^2(\nu g) \sum_{\{P \mid \nu\delta^{-1} \in P(F) \text{ et }P_1\subset P\subset P_2 \}}  \eps_P^G\cdot k_{P,X}(\delta g)|
$$
qu'on majore par
$$\sum_{\delta \in P_0(F)\back G(F)}  \sum_{P_1\subset P_2} \sum_{\nu\in P_1(F)\back G(F)}F^1(\nu g) \sigma_1^2(\nu g) |\sum_{\{P \mid \nu\delta^{-1} \in P(F) \text{ et }P_1\subset P\subset P_2 \}}  \eps_P^G\cdot k_{P,X}(\delta g)|
$$
$$
= \sum_{P_1\subset P_2} \sum_{\nu\in P_1(F)\back G(F)}F^1(\nu g) \sigma_1^2(\nu g) \sum_{\delta \in P_0(F)\back G(F)}|\sum_{\{P \mid \nu\delta^{-1} \in P(F) \text{ et }P_1\subset P\subset P_2 \}}  \eps_P^G\cdot k_{P,X}(\delta g)|
$$
qui, par le changement de variable $\delta'=\delta \nu^{-1}$, est égal à 
$$
= \sum_{P_1\subset P_2} \sum_{\nu\in P_1(F)\back G(F)}F^1(\nu g) \sigma_1^2(\nu g) \sum_{\delta \in P_0(F)\back G(F)}|\sum_{\{P \mid \delta \in P(F) \text{ et }P_1\subset P\subset P_2 \}}  \eps_P^G\cdot k_{P,X}(\delta \nu g)|
$$
Comme au §\ref{S:red1}, on  voit donc qu'il suffit de prouver  pour $P_1\subsetneq P_2$ la convergence de l'expression suivante
$$\int_{P_1(F)\back G(\AAA)^1}  F^1( g) \sigma_1^2(g) \sum_{\delta \in P_0(F)\back G(F)}|\sum_{\{P \mid \delta \in P(F) \text{ et }P_1\subset P\subset P_2 \}}  \eps_P^G\cdot k_{P,X}(\delta g)|\, dg
$$

  Soit $\al\in \Delta_1^2$ et $Q$ le sous-groupe parabolique maximal contenant $Q$ tel que $\Delta_1^Q=\Delta_1-\{\al\}$. On reprend les notations  du §\ref{S:reduc3} dans la situation $P_1=P_2$ à ceci près : les groupes notés $P_2$ et $P_3$ là-bas sont notés ici respectivement $P_1$ et $P_2$. On fixe $P\subsetneq P'$ deux sous-groupes paraboliques compris entre $P_1$ et $P_2$ qui ont même trace sur $M_Q$. Soit $\chi_P$ et $\chi_{P'}$ les fonctions caractéristiques respectives de $P(F)$ et $P'(F)$.

Par le même raisonnement que dans la preuve du théorème \ref{eq:asymp}, on voit qu'il suffit de prouver le lemme suivant
  
\begin{lemme}\label{lem:reduc-finale}
Soit $\kc$ une partie compacte de $\Sc(\ggo(\AAA))$. Il existe $c>0$ tel que pour tout $f\in \kc$ et tout $a\in A_1^{2,+}$ (cf. \eqref{eq:A13}), 

$$
 \exp(\bg -2\rho_1,H_1(a)\bd) \sum_{\delta \in P_0(F)\back G(F)}|  \chi_{P'}(\delta) \cdot k_{P',X}(\delta a)-  \chi_P(\delta) \cdot k_{P,X}(\delta a)|     \leq c\cdot \al(a)^{-1}.
$$
\end{lemme}

\begin{preuve}
On fixe une fonction $\phi\in \Sc(\ggo(\AAA))$ telle que 
$$|f|\leq \phi$$
pour tout $f\in \kc$.

Si $\delta\notin P'(F)$ ou si $P_0\not\subset P'$ l'expression en valeur absolue est nulle. Donc on peut se limiter au cas où $P_0\subset P' $ et où  $\delta\in P_0(F)\back P'(F)\simeq (P_0\cap M')(F)\back M'(F)$. On examine donc l'expression
 $$
 \exp(\bg -2\rho_1,H_1(a)\bd) \sum_{\delta \in (P_0\cap M')(F)\back M'(F)}|  k_{P',X}(\delta a)-  \chi_P(\delta) \cdot k_{P,X}(\delta a)| .    
$$

Supposons tout d'abord $P_0\not\subset P$. L'expression ci-dessus se simplifie en 
$$
 \exp(\bg -2\rho_1,H_1(a)\bd) \sum_{\delta \in  (P_0\cap M')(F)\back M'(F)}|  k_{P',X}(\delta a)|     
$$
qu'on majore par 
$$ \exp(\bg -2\rho_1^{P'},H_1(a)\bd) \sum_{X'} \int_{\ngo_{P'}(\AAA)} \phi((\Ad(a)^{-1}X)+U)\,dU
$$
où l'on somme sur les éléments $X'$ dans la classe de $M'(F)$-conjugaison de $X_{P'}$. Distinguons dans cette somme deux contributions :
\begin{enumerate}
\item celles des $X'$ qui n'appartiennent pas à $\pgo$ : on traite leur contribution par le lemme \ref{lem:maj1}.
\item celles des $X'$ qui appartiennent à $\pgo$.
\end{enumerate}
Considérons ce dernier cas : soit $X'$ un tel élément qu'on écrit $Y+V$ avec $Y\in \mgo_P(F)$ et $V\in \ngo_P^{P'}(F)$. Soit $R$ un sous-groupe parabolique standard tel que $Y\in I_R^{P}(0)$. Si $X'\in I_P^{P'}(Y)$ alors par transitivité de l'induction 
$$I_{P_0}(0)=I_{P'}(X_{P'})=I_{P'}(X')=I_{P'}(I_{P}^{P'}(Y))=I_{P}^{G}(Y)=I_{R}^G(0)$$
donc nécessairement $P_0$ et $R$ sont associés donc égaux. Donc $P_0\subset P$ ce qui n'est pas. Ainsi, l'ensemble sur lequel on somme dans cette second contribution est inclus dans l'ensemble des $Y+V$ avec $Y\in \mgo_P(F)$ et $V\in \ngo_P^{P'}(F)$ tel que $Y+V\notin I_P^{P'}(Y)$. On raisonne comme dans la preuve du lemme \ref{lem:reduc3} (cf. en particulier la majoration de \eqref{eq:p4}).

Supposons ensuite $P_0\subset P$. On traite d'abord la contribution des $\delta\notin P(F)$. Cette contribution se majore par 
$$ \exp(\bg -2\rho_1^{P'},H_1(a)\bd) \sum_{X'} \int_{\ngo_{P'}(\AAA)} \phi((\Ad(a)^{-1}X)+U)\,dU
$$
où l'on somme sur  les éléments $X'$ qui sont conjugués à $X_{P'}$ par un élément $\delta\in M'(F)-(M'\cap P)(F)$. De nouveau, on considère deux sortes de contribution : celles des éléments $X'\notin \pgo(F)$ et la contribution des éléments de $\pgo(F)$. La première se majore comme précédemment. Considérons $Y+V$ avec $Y\in \mgo_P(F)$ et $V\in \ngo_P^{P'}(F)$ qui est conjugué à $X_{P'}$ par un élément $\delta\in M'(F)-(M'\cap P)(F)$. Si $Y+V\in I_P^{P'}(Y)$ alors, comme précédemment, on voit que $Y\in I_{P_0}^P(0)$ donc que $Y$ et $X_P$ appartiennent à la même $M$-orbite notée $\of$. Par conséquent $X_{P'}$ et $Y+V$ appartiennent à $(\of\oplus \ngo_P^{P'})\cap I_P^{P'}(\of)$. Cette intersection est en fait une $M'\cap P$-orbite. Il existe $\nu\in M'\cap P$ tel que
$$\Ad(\nu)X_{P'}=Y+V=\Ad(\delta)X_{P'}$$
il s'ensuit que $\delta\in \nu M'_{X_P'}\subset PP_0=P$ ce qui est une contradiction. Donc on est ramené à une situation déjà traitée (cf. la majoration de \eqref{eq:p4}).

Il nous reste pour $P_0\subset P$ à traiter la contribution des $\delta\in P(F)$. Soit $\delta \in  P(F)$. On a 
\begin{equation}
  \label{eq:kp1}
  k_{P,X}(\delta a)=\sum_{\nu \in M_{X_P}(F)\back (M_P\cap P_0)(F)} \int_{\ngo_P(\AAA)} f(\Ad(\delta a)^{-1}(\Ad(\nu)^{-1}X_P+U)\, dU.
\end{equation}

 Comme dans la preuve du lemme \ref{lem:dec}, à l'aide d'une formule sommatoire de Poisson, on peut remplacer cette expression par 

 \begin{equation}
   \label{eq:kp2}
\tilde{k}_{P,X}(\delta a) =  \sum_{\nu \in M_{X_P}(F)\back (M_P\cap P_0)(F)}\sum_{V\in \ngo_P^{P'}(F)} \int_{\ngo_{P'}(\AAA)} f(\Ad(\delta a)^{-1}(\Ad(\nu)^{-1}X_P+V+U)\, dU,
 \end{equation}
où  $M_{X_P}$ est le centralisateur de $X_P$ dans $M_P$. Les deux expressions  \eqref{eq:kp1} et  \eqref{eq:kp2} ne sont pas égales. Cependant, la somme sur  $\delta\in P_0(F)\back P(F)$ de leurs différences satisfait la majoration cherchée (cf. lemme \ref{lem:maj2}). On doit donc majorer
\begin{equation}
  \label{eq:kp4}
  \exp(\bg -2\rho_1,H_1(a)\bd) \sum_{ \delta\in P_0(F)\back P(F)  } |k_{P',X}(\delta ) - \tilde{k}_{P,X}(\delta a) |
\end{equation}
où
\begin{equation}
  \label{eq:kp3}
  k_{P',X}(\delta a)=\sum_{\nu \in M_{X_{P'}}(F)\back (M_{P'}\cap P_0)(F)} \int_{\ngo_{P'}(\AAA)} f(\Ad(\delta a)^{-1}(\Ad(\nu)^{-1}X_{P'}+U)\, dU.
\end{equation}

 Observons que dans \eqref{eq:kp2}, on a 
$$\Ad(\nu)^{-1}X_P+V\in I_P^{P'}(X_P)
$$
si et seulement si $\Ad(\nu)^{-1}X_P+Y$ est $(M_{P'}\cap P)(F)$-conjugué à $X_{P'}$. Par conséquent, on a 
$$
 \tilde{k}_{P,X}(\delta a)-k_{P',X}(\delta a)=  \sum_{\nu \in M_{X_P}(F)\back (M_P\cap P_0)(F)}$$
 $$\sum_{\{V\in \ngo_P^{P'}(F)\mid \Ad(\nu)^{-1}X_P+V\notin I_P^{P'}(X_P)} \int_{\ngo_{P'}(\AAA)} f(\Ad(\delta a)^{-1}(\Ad(\nu)^{-1}X_P+V+U)\, dU
$$
Finalement, l'expression \eqref{eq:kp4} se majore par
$$
 \exp(\bg -2\rho_1,H_1(a)\bd)\sum_Y \sum_{\{V\in \ngo_P^{P'}(F)\mid Y+V\notin I_P^{P'}(Y)} \int_{\ngo_{P'}(\AAA)} f(\Ad(a)^{-1}(Y+V+U)\, dU
$$
et l'on conclut comme dans la preuve du lemme \ref{lem:reduc3}.
\end{preuve}

\end{paragr}

\section{Calcul de coefficients globaux d'Arthur.}\label{sec:coef}

\begin{paragr} Les notations sont celles des sections \ref{sec:IOP} et  \ref{sec:reg}. Soit $S$ un ensemble fini de places. Arthur a défini dans \cite{wei_or} des intégrales orbitales pondérées locales unipotentes. Sa construction s'adapte au cas des algèbres de Lie. \emph{A priori}, ces distributions apparaissent comme des limites d'intégrales orbitales pondérées plus explicites. En fait, et c'est le point de vue adopté dans \cite{scuft}, on peut extraire de \cite{wei_or} une expression intégrale pour ces distributions. Rappelons brièvement leur définition. Soit $f\in \Sc(\ggo(\AAA_S))$, $L\in \lc(T)$ et $\of$ une $L$-orbite nilpotente dans $\lgo$ : on définit
$$
J_L^G(\of,f)=\int_{\ngo_P(\AAA_S)} \int_{K_S} f(\Ad(k^{-1})U) w_L(U) \,dU.
$$
où $P$ est un certain sous-groupe parabolique et $w_L$ une fonction définie à l'aide de la théorie des $(G,M)$-familles. Pour des définitions détaillées, on renvoie à  \cite{scuft} §§ 8.5 et 12.5. Précisons cependant que l'intégrale dépend de  choix de mesure : on suit ceux adoptés dans le présent article qui  ne sont pas exactement ceux de \cite{scuft}.
\end{paragr}

\begin{paragr}[Développement d'Arthur.] --- Il s'agit de développer la distribution globale $J_\of(f)$ en termes des distributions nilpotentes pondérées introduites ci-dessus.
  
  \begin{theoreme}\label{thm:coef}
    Soit $\of$ l'orbite de $X$ et $S$ un ensemble fini de places. Pour toute fonction $f_S\in \Sc(\ggo(\AAA_S))$, on a 
    \begin{equation}
      \label{eq:JOX}
       J_{\of}(f_S\otimes \mathbf{1}^S)=\vol([M_0]^1)\cdot \sum_{L\in \lc(M_0)}\frac{|W_0^L|}{|W_0|} a^L(S,\of_L) J_L^G(\of_L, f_S)
    \end{equation}
  \end{theoreme}
où 
\begin{enumerate}
\item $\of_L=\Ind_{M_0}^L(0)$ ;
\item $ a^L(S,\of_L)$ est la valeur en $0$ de l'expression suivante (bien définie pour  $\la$  dans un ouvert de $a_{M_0,\CC}^{L,*}$)
 $$\frac{1}{|W_0^L|}\sum_{w\in W_0^L} \varphi(w\cdot \la) \theta_{0}^L(w\la)
$$
où $ \theta_{0}^L(\la)=\theta_{P_0\cap}^L$ est défini au §\ref{S:theta} et où
$$\varphi(\la)=  \prod_{\al\in \Delta_{P_0\cap L}^L} \frac{ \tilde{Z}_d^S(d+\frac1d \bg \la,\varpi_\al^\vee\bd)  }{ \tilde{Z}_d^S(d) }
$$
pour $\la\in a_{M_0,\CC}^{L,*}$.
\end{enumerate}

\begin{remarque} En exploitant la proposition \ref{prop:cp} et son corollaire \ref{cor:cp}, on peut donner deux expressions alternatives pour le coefficient  $ a^L(S,\of_L)$ : c'est la valeur en $0$ qui est commune aux deux expressions suivantes 
$$\sum_{ Q} \eps_Q^L \cdot \hat{\theta}_{P_0\cap L}^Q(\la)^{-1} \varphi(\la^Q) \theta_Q(\la)^{-1}
$$
et
$$\sum_{ Q}  \eps_{P_0}^P \cdot \hat{\theta}_{P_0\cap L}^P(\la)^{-1} \varphi(\la_Q) \theta_Q(\la)^{-1}
$$
  où l'on somme sur les sous-groupes paraboliques $Q$ de $L$ qui contiennent $P_0\cap L$.
\end{remarque}

\end{paragr}

\begin{paragr}[Variante.] --- \label{S:variante} Les distributions $J_L^G(\of_L)$ ne dépendent que la $W$-orbite de $(L,\of_L)$. On voit qu'il en est de même des coefficients $a^{L}(\of,S)$. Soit $(L,\of_L)$ tel que $I_L^G(\of_L)=\of$. Nécessairement, il existe $w\in W$ tel que $wM_0w^{-1}\subset L$ et $\of_L=I_{wM_0w^{-1}}^L(0)$. On pose alors
$$\tilde{a}^{L}(S,\of_L)= \vol([M_0]^1)\cdot a^{w^{-1}Lw}(S, \of_{ w^{-1}Lw }).$$
Il vient alors 
\begin{eqnarray*}
   J_{\of}(f_S\otimes \mathbf{1}^S)&=& \sum_{L\in \lc(M_0)/W_0} |W_0(L)|^{-1} \tilde{a}^L(S,\of_L) J_L^G(\of_L, f_S)\\
&=&  \sum_{L\in \lc(T)/W} |W(L)|^{-1} \sum_{\of' }  \tilde{a}^L(S,\of') J_L^G(\of', f_S)\\
&=&  \sum_{L\in \lc(T)}\frac{|W^L|}{|W|} \sum_{\of' }  \tilde{a}^L(S,\of') J_L^G(\of', f_S)\\
 \end{eqnarray*}
  où
  \begin{itemize}
  \item on somme sur les $L$-orbites nilpotentes $\of'$ dans $\lgo$ tels que $I_L^G(\of')=\of$ ;
  \item on pose $W(L)=\Norm_{W}(L)/W^L$ et  $W_0(L)=\Norm_{W_0}(L)/W_0^L$ si de plus $M_0\subset L$ ; dans ce dernier cas, ces groupes sont isomorphes. 
  \end{itemize}
On retrouve formellement le développement d'Arthur (corollaire 8.4 de \cite{ar_unipvar}). Au passage à l'algèbre de Lie près, est-ce le même développement terme à terme ? Par construction, les distributions  $J_L^G(\of')$ sont les mêmes. On peut montrer que les coefficients  $\tilde{a}^L(S,\of')$ sont aussi ceux définis par Arthur (cf. la remarque 12.7.2 de \cite{scuft}). En particulier, pour $L=M_0$, la seule orbite $\of'$ possible est l'orbite nulle $(0)$ et on trouve
$$   \tilde{a}^{M_0}(S,(0)) =\vol([M_0]^1)$$
comme calculé par Arthur (corollaire 8.5 de  \cite{ar_unipvar}).

\end{paragr}

\begin{paragr}[Preuve du théorème \ref{thm:coef}.] L'expression obtenue au théorème \ref{ref:calcul} repose sur la combinatoire d'une $(G,M_0)$-famille d'Arthur obtenue par une construction simple à partir de la fonction $J_{G,X}(f,\la)$ ou plutôt de sa variante $\tilde{J}_{G,X}(f,\la)$ (cf. § \ref{S:cstr}). On a pour $\la$ dans un ouvert de $a_{P_0,\CC}^{G,*}$ et $f=f_S\otimes\mathbf{1}_\ggo^S$
  \begin{eqnarray*}
    \tilde{J}_{G,X}(f,\la) &=& \vol([G_X]^1)\hat{\theta}^G_0(\la) \int_{G_X(\AAA)\back G(\AAA)} \exp(-\bg\la,H_0(g)\bd) f(\Ad(g^{-1})(X)) \, dg\\
&=&  \vol([G_X]^1)d^{|\Delta_0^G|}\vol(a_0^G/\ZZ(\hat{\Delta}_0^{G,\vee}))^{-1} \cdot \psi(\la) \phi(\la)
  \end{eqnarray*} 
où l'on pose
 $$\psi(\la) =  \int_{G_X(\AAA_S)\back G(\AAA_S)} \exp(-\bg\la,H_0(g)\bd) f_S(\Ad(g^{-1})(X)) \, dg$$
et
\begin{eqnarray*}
  \phi(\la)&=& d^{-|\Delta_0^G|}\vol(a_0^G/\ZZ(\hat{\Delta}_0^{G,\vee}))\hat{\theta}^G_0(\la) \int_{G_X(\AAA^S)\back G(\AAA^S)} \exp(-\bg\la,H_0(g)\bd) \mathbf{1}^S_\ggo(\Ad(g^{-1})(X)) \, dg\\
&=&  \prod_{\al\in \Delta_0^G} \tilde{Z}_d^S(d+\frac1d \bg \la,\varpi_\al^\vee\bd).
\end{eqnarray*}
La seconde égalité est un calcul analogue à celui de la proposition \ref{prop:cal-unite} assertion 2. On introduit alors les $(G,M_0)$-familles définies pour $w\in W_0$ par 
$$\psi_{wP_0}=\psi(w^{-1}\la)   \text{  et  } \phi_{wP_0}=\phi(w^{-1}\la).
$$
On a d'ailleurs pour $P\in \pc(M_0)$
$$
\phi_P(\la)=\prod_{\al\in \Delta_{P}^G} \tilde{Z}_d^S(d+\frac1d \bg \la,\varpi_\al^\vee\bd)
$$

Soit $L\in \lc(M_0)$. On sait associer à $\phi$ et  $Q\in \pc(L)$ une $(L,M)$-famille : pour tout parabolique $R\in \pc^{L}(M_0)$ et $\la\in a_{M_0,\CC}^{L,*} $ on pose
$$\phi_R^{Q}(\la)=\varphi_{RN_{Q}}(\la).
$$
On a donc explicitement 
$$\phi_R^{Q}(\la)= \tilde{Z}_d^S(d)^{|\Delta_L^G|}\prod_{\al\in \Delta_{R}^L} \tilde{Z}_d^S(d+\frac1d \bg \la,\varpi_\al^\vee\bd).
$$
En particulier, on voit que $\phi_R^{Q}(\la)$ ne dépend pas du choix de $Q\in \pc(L)$. On note simplement $\phi_R^L$ cette fonction. De même, à $\psi$ on associe une $(G,L)$-famille : pour $\la\in a_{L,\CC}^{G,*}$ et $Q\in \pc^G(L)$ on pose
$$\psi_{Q}(\la)=\psi_P(\la)$$
pour un quelconque élément $P\in \pc^{Q}(M_0)$ (le résultat ne dépend pas du choix de $P$). Finalement, on note $\phi_{M_0}^L$ et $\psi_L$ les valeurs respectives en $\la=0$ des fonctions (lisses) suivantes
$$
\la\in  a_{M_0,\CC}^{L,*}\mapsto  \sum_{P\in \pc^L(M_0)}  \phi_P^L(\la) \theta_P^L(\la)^{-1}
$$
et
$$
\la\in  a_{L,\CC}^{G,*}\mapsto  \sum_{P\in \pc^G(L)}  \psi_P(\la)) \theta_P^G(\la)^{-1}.
$$

Au facteur $\frac{1}{|W_0|}  \vol([G_X]^1)d^{|\Delta_0^G|}\vol(a_0^G/\ZZ(\hat{\Delta}_0^{G,\vee}))^{-1}$ près, la valeur de $J_\of(f)$ est donnée par la valeur  notée $(\psi \phi)_{M_0}$ en $0$ de l'expression 
$$\sum_{P\in \pc(M_0)}  \phi_P(\la)\psi_P(\la) \theta_P(\la)^{-1}.
$$
D'après Arthur (cf. \cite{trace_inv} corollaire 6.5), on a la formule de décomposition
$$
 (\psi \phi)_{M_0}=\sum_{L \in \lc(M_0)}  \phi_{M_0}^{L} \psi_L.
$$

D'après le lemme \ref{lem:IOPloc} ci-dessous, on a 
 $$\psi_L=  Z_{d,S}(d)^{|\Delta_0^G|} J_L^G(\of_L,f_S).$$
D'autre part, $ \phi_{M_0}^{L}$ est la valeur en le point $0$ de la fonction $\phi_{M_0}^{L}(\la)$ de la variable $\la\in a_{M_0,\CC}^{L,*}$ définie ainsi
\begin{eqnarray*}
   \phi_{M_0}^{L}& =&  \tilde{Z}_d^S(d)^{|\Delta_0^G|} \sum_{P\in \pc^L(M_0)}\prod_{\al\in \Delta_{P}^L} \frac{ \tilde{Z}_d^S(d+\frac1d \bg \la,\varpi_\al^\vee\bd)  }{ \tilde{Z}_d^S(d) }  \theta_P^L(\la)^{-1}\\
   &=&\tilde{Z}_d^S(d)^{|\Delta_0^G|}  \sum_{w\in W_0^L}  \prod_{\al\in \Delta_{P_0\cap L}^L} \frac{ \tilde{Z}_d^S(d+\frac1d \bg w^{-1}\cdot \la,\varpi_\al^\vee\bd)  }{ \tilde{Z}_d^S(d) }  \theta_{P_0\cap L}^L(w^{-1}\cdot\la)^{-1}\\
\end{eqnarray*}
Le théorème \ref{thm:coef} s'en déduit compte tenu de l'égalité
$$
 \vol([G_X]^1)d^{|\Delta_0^G|}\vol(a_0^G/\ZZ(\hat{\Delta}_0^{G,\vee}))^{-1}\cdot \tilde{Z}_d(d)^{|\Delta_0^G|}=\vol([M_0]^1)
$$
qui est donnée par le lemme \ref{lem:volGvolM} pour $P=P_0$.

\begin{lemme}\label{lem:IOPloc}
  Avec les notations ci-dessus, on a 
$$\psi_{L}= Z_{d,S}(d)^{|\Delta_0^G|} J_L^G(\of_L,f_S).
$$
\end{lemme}

\begin{preuve}
On va d'abord faire le lien avec d'autres intégrales introduites dans \cite{scuft}. Soit $P\in \pc(M_0)$ et $w\in W_0$ tel que $w^{-1}Pw =P_0$. Pour tout $\la \in a_{0,\CC}^{G,*}$ et $g\in G(\AAA_S)$, on a 
$$\bg w^{-1}\cdot \la, H_{P_0}(g)\bd=\bg  \la, w\cdot H_{P_0}(g)\bd = \bg  \la, H_{P_0}(w g)\bd $$
(on profite ici du fait que $w$ peut être représenté par une matrice de permutation); dans la terminologie de \cite{scuft}, le groupe $P_0$ est un groupe de Richardson pour l'élément $X$. Le vecteur  $H_{P_0}(w g)$ y est noté $R_P(g)$ (\emph{ibid.} § 5.3). Plus généralement, pour tout $Q\in \fc(M_0)$, on peut définir un vecteur $R_Q(g)\in a_{M_Q}$ (\emph{ibid.} § 5.4). On vérifie, et c'est immédiat sur les définitions, qu'on a pour tout  $\la$ dans un ouvert de $a_{L,\CC}^{G,*}$ et $Q\in \pc(L)$ 
 $$\psi_Q(\la) =  \int_{G_X(\AAA_S)\back G(\AAA_S)} \exp(-\bg\la,R_Q(g)\bd) f_S(\Ad(g^{-1})(X)) \, dg.
$$
Le membre de droite fait sens ici car pour $h\in  G_X(\AAA_S)$ on a $R_Q(hg)=R_Q(h)+R_Q(g)$ (\emph{ibid.} lemme 5.6.1) et $R_Q(h)\in a_G$ ici.

  Soit $k=\dim(a_L^G)$ et $\la\in a_{L,\CC}^{G,*}$ un point assez général. Il résulte de la théorie des $(G,M)$-familles (cf. \cite{trace_inv} formule (6.5) p. 37) que 
$$\psi_{L}=  \frac{1}{k!} \sum_{Q\in \pc^G(L)}  \frac{d^k}{dt^k}(\psi_Q(t \la))_{|t=0}\, \theta_Q^G(\la)^{-1}.
$$
En dérivant sous l'intégrale, il vient
$$\frac{d^k}{dt^k}(\psi_P(t \la))= \int_{G_X(\AAA_S)\back G(\AAA_S)} (-\bg\la,R_Q(g)\bd)^k \exp(-\bg t \la,R_Q(g)\bd) f_S(\Ad(g^{-1})(X)) \, dg
$$
et on peut évaluer cette expression en $t=0$ (cf. les énoncés de convergence de la section 8 de \emph{ibid.}).
Par conséquent, on a 
$$\psi_{L}= \int_{G_X(\AAA_S)\back G(\AAA_S)}  f_S(\Ad(g^{-1})(X)) v_{L,X}^G(g) \, dg
$$
où le poids  est donné par la formule 
$$
v_{L,X}^G(g)= \frac{1}{k!} \sum_{Q\in \pc^G(L)}   (-\bg\la,R_Q(g)\bd)^k \exp(-\bg\la,R_Q(g)\bd)  \theta_Q^G(\la)^{-1}.
$$
Le poids est donc celui défini aux paragraphes 8.3 et 12.5 de   \emph{ibid.}. Il résulte du théorème 8.5.1 de  \emph{ibid.} qu'avec \emph{nos} choix de mesures on a 
$$
\psi_{L}= Z_{d,S}(d)^{|\Delta_0^G|}\cdot  J_L^G(\of_L,f_S).
$$
\end{preuve}

\end{paragr}

\section{Un résultat sur les $(G,M)$-familles}\label{sec:GM}

\begin{paragr}
  La situation dans ce paragraphe est tout-à-fait générale (c'est celle de la section \ref{sec:notations}). Soit $P_0\subset G$ un sous-groupe parabolique.

  \begin{proposition}\label{prop:cp}  Soit $\varphi$ une fonction lisse sur $a_{P_0}^{G,*}$. Les fonctions     
$$\tilde{c}_{P_0}(\varphi,\la)=\sum_{P_0\subset P\subset G} \eps_P^G \cdot \hat{\theta}_0^P(\la)^{-1} \varphi(\la^P) \theta_P(\la)^{-1}
$$
et
$$c_{P_0}(\varphi,\la)=\sum_{P_0\subset P\subset G} \eps_{P_0}^P \cdot \hat{\theta}_0^P(\la)^{-1} \varphi(\la_P) \theta_P(\la)^{-1}
$$
sont lisses en $\la=0$.
En outre,
\begin{equation}
  \label{eq:cp}
  c_{P_0}(\varphi,0)=\tilde{c}_{P_0}(\varphi,0).
\end{equation}
  \end{proposition}

  \begin{preuve}
Soit $k=\dim(a_{P_0}^G)$. La lissité de la fonction $c_{P_0}$ en $0$ est un résultat d'Arthur (cf. lemme 6.1 de \cite{trace_inv}). La lissité de $\tilde{c}_{P_0}$ en $0$ s'obtient par une variante des méthodes d'Arthur (cf. la preuve du théorème 2.6.4 de \cite{scfhn}). En utilisant un développement de Taylor, on voit que pour tout point $\la$ assez général, on a
$$\tilde{c}_{P_0}(\varphi,0)=\frac{1}{k!}\sum_{P_0\subset P\subset G} \eps_P^G \cdot \hat{\theta}_0^P(\la)^{-1}  (\frac{d^k}{dt^k}\varphi(t \la^P))_{|t=0} \theta_P(\la)^{-1}
$$
et
$$c_{P_0}(\varphi,0)=\frac{1}{k!}\sum_{P_0\subset P\subset G}  \eps_{P_0}^P \cdot \hat{\theta}_0^P(\la)^{-1} (\frac{d^k}{dt^k}\varphi(t \la_P))_{|t=0} \theta_P(\la)^{-1}.
$$
En particulier, chaque expression ne dépend que du comportement local de $\varphi$ en $0$. On peut donc et on va supposer que $\varphi$ est lisse et à support compact. On peut l'écrire alors comme la transformée de Fourier d'une fonction $\hat{\varphi}$ lisse et à décroissance rapide sur $a_{P_0}^{G}$ :
$$
\varphi(\la)=\int_{\RR}  \hat{\varphi}(H) \exp(i \bg \la,H\bd ) \, dH
$$
On a donc 
$$ (\frac{d^k}{dt^k}\varphi(t \la))_{|t=0}=\int_{\RR}  \hat{\varphi}(H)    (\frac{d^k}{dt^k}\exp(i t\bg \la,H\bd )_{|t=0} \, dH.
$$
Il vient donc
$$\tilde{c}_{P_0}(\varphi,0)= \int_{\RR}  \hat{\varphi}(H) \cdot \tilde{c}_{P_0}(  \exp(i \bg \la,H\bd ),0) \, dH
$$
et
$$c_{P_0}(\varphi,0)= \int_{\RR}  \hat{\varphi}(H) \cdot c_{P_0}(  \exp(i \bg \la,H\bd ),0) \, dH.
$$

Il suffit donc de prouver \eqref{eq:cp} pour les fonctions $\varphi$ telles que $\varphi(\la)=\exp(i \bg \la,H\bd )$. Pour une telle fonction, on a (en tenant compte de  $\eps_{P_0}^G=(-1)^k$)
\begin{eqnarray*}
   \tilde{c}_{P_0}(  \exp(i \bg \la,H\bd ),\la)&=&  \eps_{P_0}^G \exp(i \bg \la,H\bd ) \cdot c_{P_0}(  \exp(-i \bg \la,H\bd ),\la)\\
&=& \exp(i \bg \la,H\bd ) \cdot c_{P_0}(  \exp(i \bg \la,H\bd ),- \la)\\
\end{eqnarray*}
En $\la=0$, on obtient l'égalité \eqref{eq:cp} pour la fonction  $\exp(i \bg \la,H\bd )$ : cela conclut.
  \end{preuve}

\end{paragr}

\begin{paragr}  Plaçons-nous dans le cas de $G=GL(n)$ et $P_0$ le sous-groupe parabolique défini au§\ref{S:P0}. Soit $P_0=M_0N_0$ la décomposition de Levi standard  de $P_0$. On peut  formuler le corollaire suivant à la proposition \ref{prop:cp}.

\begin{corollaire}\label{cor:cp}Soit $\varphi$ une fonction sur $a_{P_0,\CC}^{G,*}$ holomorphe au voisinage de $0$. L'expression
  \begin{equation}
    \label{eq:cM}
    \frac{1}{|W_0|}\sum_{w\in W_0} \varphi(w^{-1}\cdot \la) \theta_{0}(w^{-1}\la)
  \end{equation}
  est holomorphe en $\la=0$ et sa valeur en ce point est égale à 
$$\tilde{c}_{P_0}(\varphi,0).
$$
  
\end{corollaire}

\begin{preuve}
   Introduisons pour tout $w\in W_0$
$$\varphi_{wP_0}(\la)=\varphi(w^{-1}\cdot\la).
$$   
La famille $(\varphi_P)_{P\in \pc(M_0)}$ est alors dans la terminologie de \cite{trace_inv} une $(G,M_0)$-famille. Que l'expression \eqref{eq:cM} soit holomorphe en $\la=0$ résulte de la théorie d'Arthur (cf. preuve du lemme 6.2 de \cite{trace_inv}) et du fait qu'ici l'ensemble $\pc(M_0)$ des sous-groupes paraboliques qui admettent $M_0$ comme facteur de Levi est  un $W_0$-torseur. En outre,  par le  corollaire 6.4 de \cite{trace_inv} on a 
$$\sum_{w\in W_0} \varphi(w^{-1}\cdot \la) \theta_{0}(w^{-1}\la)= \sum_{w\in W_0  } c_{wP_0}(\varphi(w^{-1}\cdot\la), \la)
$$  
On vérifie aisément qu'on a 
 $$c_{wP_0}(\varphi(w^{-1}\cdot\la), \la)= c_{P_0}(\varphi,w^{-1}\la).
$$
Il s'ensuit que la valeur en $\la=0$ de \eqref{eq:cM} est égale à $|W_0|\cdot c_{P_0}(\varphi,0)$. On conclut par la proposition \ref{prop:cp}.
\end{preuve}
\end{paragr}

\bibliography{bibiog}
\bibliographystyle{plain}

\begin{flushleft}
Pierre-Henri Chaudouard \\
Université Paris Diderot (Paris 7) et Institut Universitaire de France\\
 Institut de Mathématiques de Jussieu-Paris Rive Gauche \\
 UMR 7586 \\
 Bâtiment Sophie Germain \\
 Case 7012 \\
 F-75205 PARIS Cedex 13 \\
 France
\medskip

Adresse électronique :\\
Pierre-Henri.Chaudouard@imj-prg.fr \\
\end{flushleft}
\end{document}